\documentclass[a4paper,12pt,reqno]{amsart}
\usepackage[dvips]{graphics} 
\newtheorem{lem}{Lemma}
\newtheorem{theo}{Theorem}
\newtheorem{defi}{Definition}
\newtheorem{rem}{Remark}
\newtheorem{cor}{Corollary}

\begin{document}

\vspace*{4cm}

\begin{center}
\LARGE \sc
A characterisation of Pfaffian near bipartite graphs
\end{center}

\vspace{2cm}

\begin{center}
\Large Ilse Fischer \\
        Universit\"at Klagenfurt \\
        A-9020 Klagenfurt\\
        Austria\\
        Email: \tt Ilse.Fischer@uni-klu.ac.at
\end{center}

\vspace{1cm}

\begin{center}
\Large  C.H.C. Little\footnote{The second author thanks the University of Klagenfurt for its 
        hospitality while this research was undertaken.} \\
        Massey University \\
        Palmerston North\\
        New Zealand \\
        Email: \tt c.little@massey.ac.nz
\end{center}

\newpage

\vspace*{3cm}

{\Large \sc  Pfaffian near bipartite graphs }

\vspace{2cm}

{\parindent0cm  Ilse Fischer \\
Universit\"at Klagenfurt \\
Universit\"atsstrasse 65-67 \\
A-9020 Klagenfurt\\
Austria\\
Email: \tt Ilse.Fischer@uni-klu.ac.at}
\newpage

\vspace*{3cm}

\begin{center}
\large Abstract
\end{center}

\vspace{1cm}
A graph is {\it 1-extendible} if every edge has a 1-factor
containing it. A 1-extendible non-bipartite graph $G$ is said to be
{\it near bipartite} if there exist edges $e_1$ and $e_2$ such that
$G - \{e_1, e_2\}$ is 1-extendible and bipartite. We characterise
the Pfaffian near bipartite graphs in terms of forbidden subgraphs.
The theorem extends an earlier characterisation of Pfaffian bipartite
graphs.

\newpage

\section{Introduction}
The graphs considered in this paper are finite and have no loops or multiple edges. They are
also undirected and connected unless an indication to the contrary is given. 
If $v$ and $w$ are vertices in a directed graph, then $(v,w)$ denotes an edge joining $v$ and 
$w$ and directed from $v$ to $w$. If $G$ is any graph, then we denote its vertex set by $VG$ and 
its edge set by $EG$. A {\it $1$-factor} of $G$ is a subset $f$ of $EG$ such that every vertex has a unique edge
of $f$ incident on it.

Let $G^*$ be a directed graph with an even number $2n$ of vertices and let $F$ be the set 
$\{f_1,f_2,\dots,f_k\}$ of $1$-factors of $G^*$. For all $i$ write
$$
f_i=\{(u_{i 1},w_{i 1}),(u_{i 2},w_{i 2}),\dots,(u_{i n},w_{i n})\},
$$
where $u_{ij},w_{ij} \in VG^*$ for all $j$. Associate with $f_i$ a plus 
sign if 
$$
u_{i 1} w_{i1} u_{i 2} w_{i 2} \dots u_{i n} w_{i n}
$$
is an even permutation of 
$$
u_{1 1} w_{1 1} u_{1 2} w_{1 2} \dots u_{1 n} w_{1 n},
$$
and a minus sign otherwise. Note that the signs of the $1$-factors are independent 
of the order in which their edges have been written. They are dependent 
on the choice of $f_1$, but the resulting partition of $F$ into 
two complementary subsets is not. If $G$ is an undirected graph, we say that $G$ is 
a {\it Pfaffian } graph if there exists an orientation such that all the $1$-factors of 
$G$ have the same sign.  
We say that this orientation is a {\it Pfaffian} orientation of $G$. Pfaffian orientations have been 
used by Kasteleyn \cite{Ka} to enumerate $1$-factors in planar graphs.
In fact his method can be used precisely for those graphs that are Pfaffian. It is therefore of interest 
to know which graphs are Pfaffian, but this question is open.

Pfaffian bipartite graphs have been characterised by Little \cite{Li1}, who proved the following
theorem. 
\begin{theo}
\label{Bipartite}
A bipartite graph $G$ is non-Pfaffian if and only if it contains an even subdivision $J$ 
of $K_{3,3}$ such that $G - VJ$ has a $1$-factor.
\end{theo}

Here we need to explain the term `even subdivision'. An {\it edge subdivision} of 
a graph $G$ is defined as a graph obtained from $G$ by replacing an edge joining vertices 
$v$ and $w$ with a path $P$ joining $v$ and $w$ but having no other vertices in common with $G$. 
The edge subdivision is {\it even} if $P$ has odd length. A graph $H$ is a 
{\it subdivision} of $G$ if for some positive integer $k$ there exist graphs $G_0,G_1,\dots,G_k$ 
such that $G_0=G$, $G_k=H$ and, for all $i>0$, $G_i$ is an edge subdivision of $G_{i-1}$. If $G_1,G_2,\dots,G_k$ 
can be chosen so that in addition $G_i$ is an even edge subdivision of $G_{i-1}$ for all $i >0$, then $H$ is said to be 
an {\it even subdivision} of $G$. It is easy to see that $G$ is Pfaffian if and only if $H$ is Pfaffian. A more general result
is proved in Lemma~\ref{degree2}. 

At this point it is worth mentioning that Robertson, Seymour and 
Thomas \cite{RoSeTh} have recently found a polynomial-time algorithm which decides 
whether a bipartite graph is Pfaffian or not.

A graph is {\it 1-extendible} if every edge has a 1-factor containing it. Such graphs are the only graphs of interest
in the study of the Pfaffian property, as any edge belonging to no 1-factor is irrelevant. 
A 1-extendible non-bipartite graph $G$ is said to be {\it near bipartite} if there exist edges $e_1$, $e_2$ such that $G - \{e_1,e_2\}$ is 
1-extendible and bipartite. If $G$ were a 1-extendible 
graph and $G - \{e\}$ were bipartite for some edge $e$, then $G$ would also be bipartite. This
observation explains why we remove two edges from $G$, rather than one, in the definition of a near bipartite graph.
The aim of the present paper is to extend Theorem \ref{Bipartite} to a characterisation of 
Pfaffian near bipartite graphs in terms of forbidden subgraphs.    

\begin{figure}
\begin{flushleft}
\setlength{\unitlength}{1cm}
\hspace{1.5cm}\mbox{\scalebox{0.35}{%
 \includegraphics{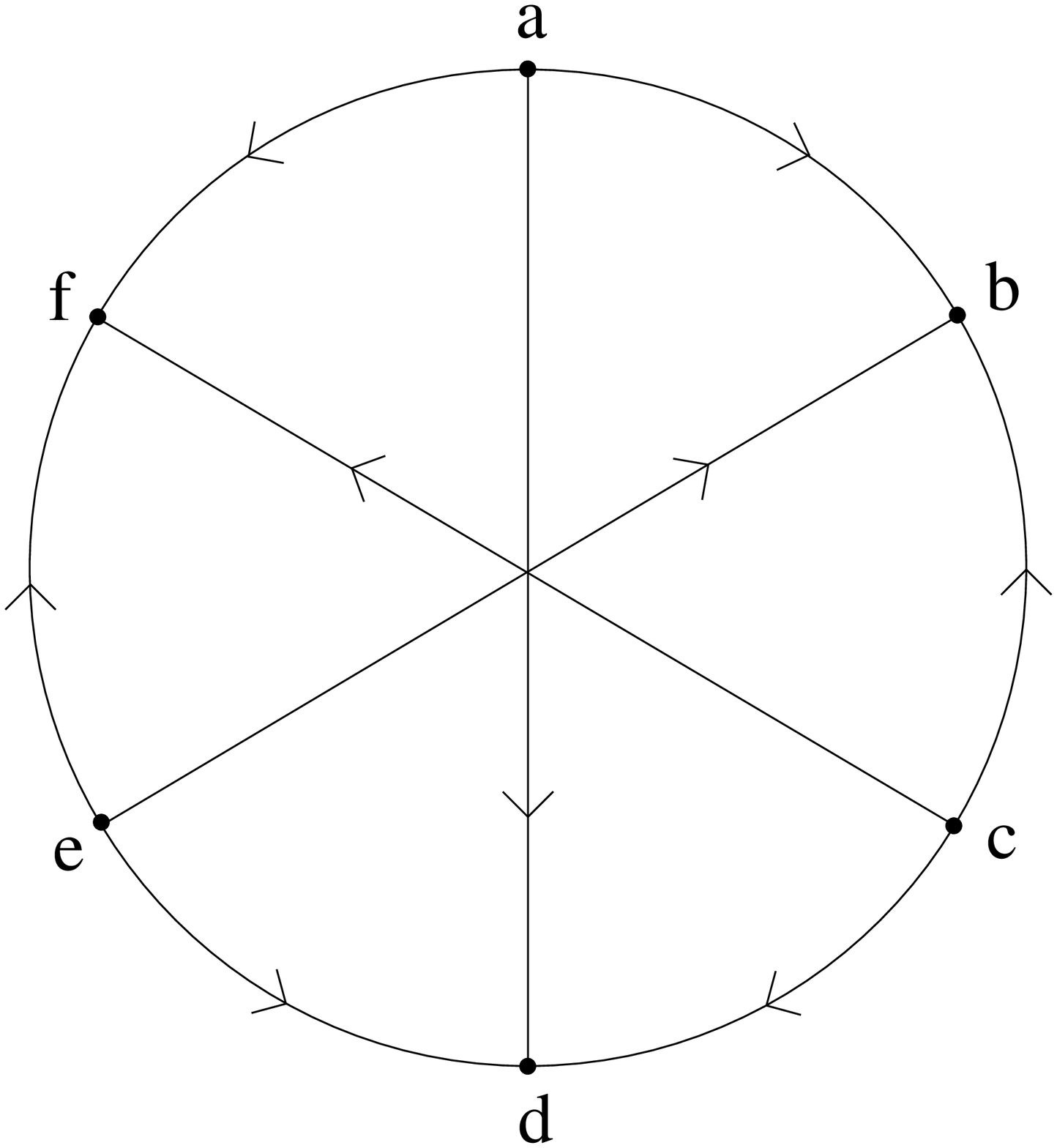}}}%
\hspace{1.5cm} \parbox[b]{1cm}{
\begin{tabular}[b]{cc}
$(a,b)(c,d)(e,f)$&$ +$\\
$(a,b)(c,f)(e,d)$&$ -$\\
$(a,d)(c,f)(e,b)$&$ +$\\
$(a,d)(c,b)(e,f)$&$ -$\\
$(a,f)(c,b)(e,d)$&$ +$\\
$(a,f)(c,d)(e,b)$&$ -$
\end{tabular}
\vspace{1.5cm}}
\end{flushleft}
\caption{The graph $K_{3,3}$.}
\label{k33}
\end{figure}

\begin{figure}
\begin{flushleft}
\setlength{\unitlength}{1cm}
\hspace{0.5cm}\mbox{\scalebox{0.35}{%
 \includegraphics{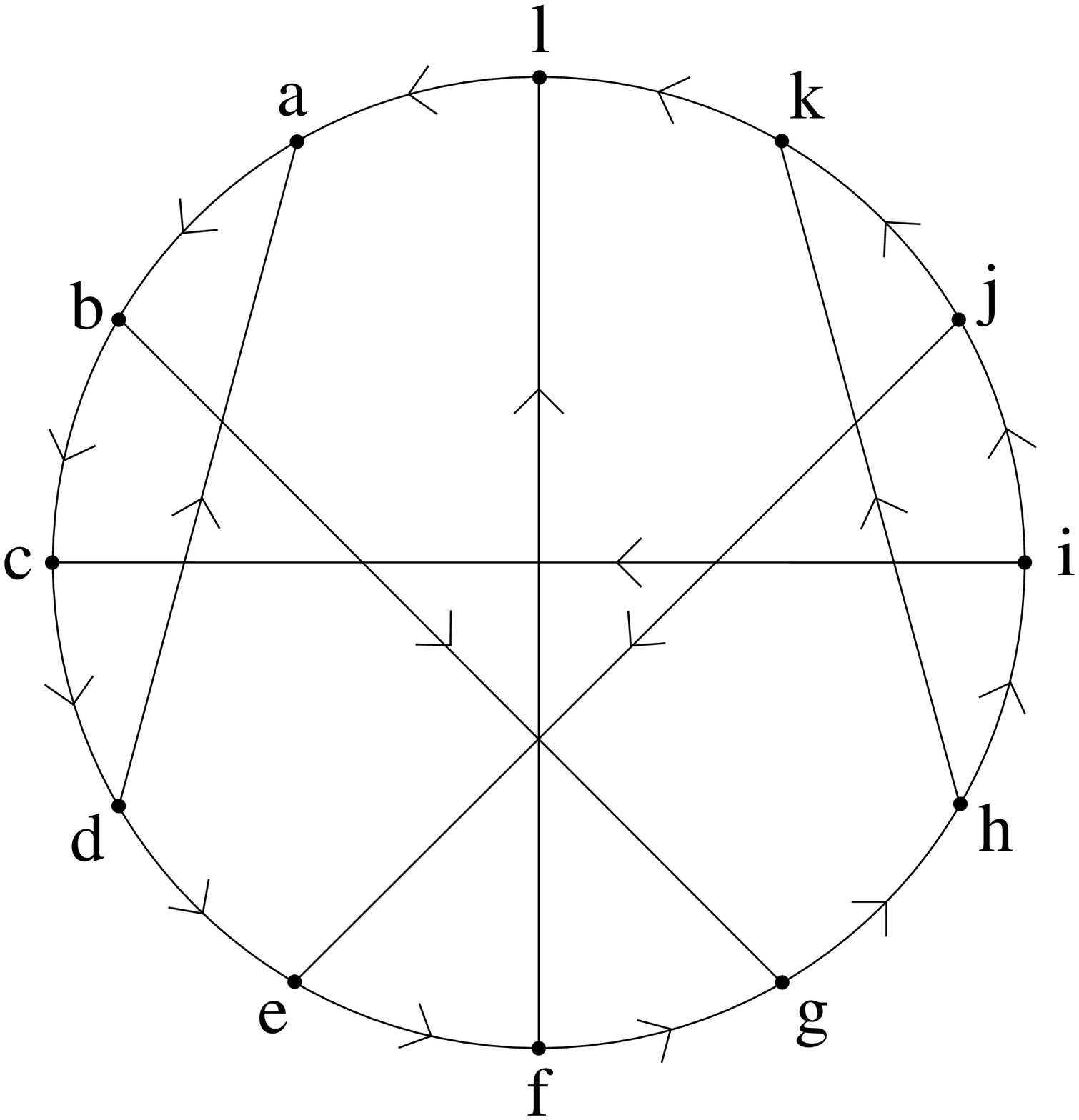}}}%
\hspace{1.5cm} \parbox[b]{1cm}{
\begin{tabular}[b]{cc}
$(a,b)(c,d)(e,f)(g,h)(i,j)(k,l)$&$ +$\\
$(b,c)(d,e)(f,g)(h,i)(j,k)(l,a)$&$ -$\\
$(d,a)(j,e)(b,c)(f,g)(h,i)(k,l)$&$ +$\\
$(b,g)(h,k)(c,d)(e,f)(i,j)(l,a)$&$ -$\\
$(i,c)(f,l)(a,b)(d,e)(g,h)(j,k)$&$ +$\\
$(d,a)(b,g)(i,c)(j,e)(h,k)(f,l)$&$ -$
\end{tabular}
\vspace{1.5cm}}
\end{flushleft}
\caption{The graph $\Gamma_1$.}
\label{gamma1}
\end{figure}

\begin{figure}
\begin{flushleft}
\setlength{\unitlength}{1cm}
\hspace{0.5cm}\mbox{\scalebox{0.35}{%
 \includegraphics{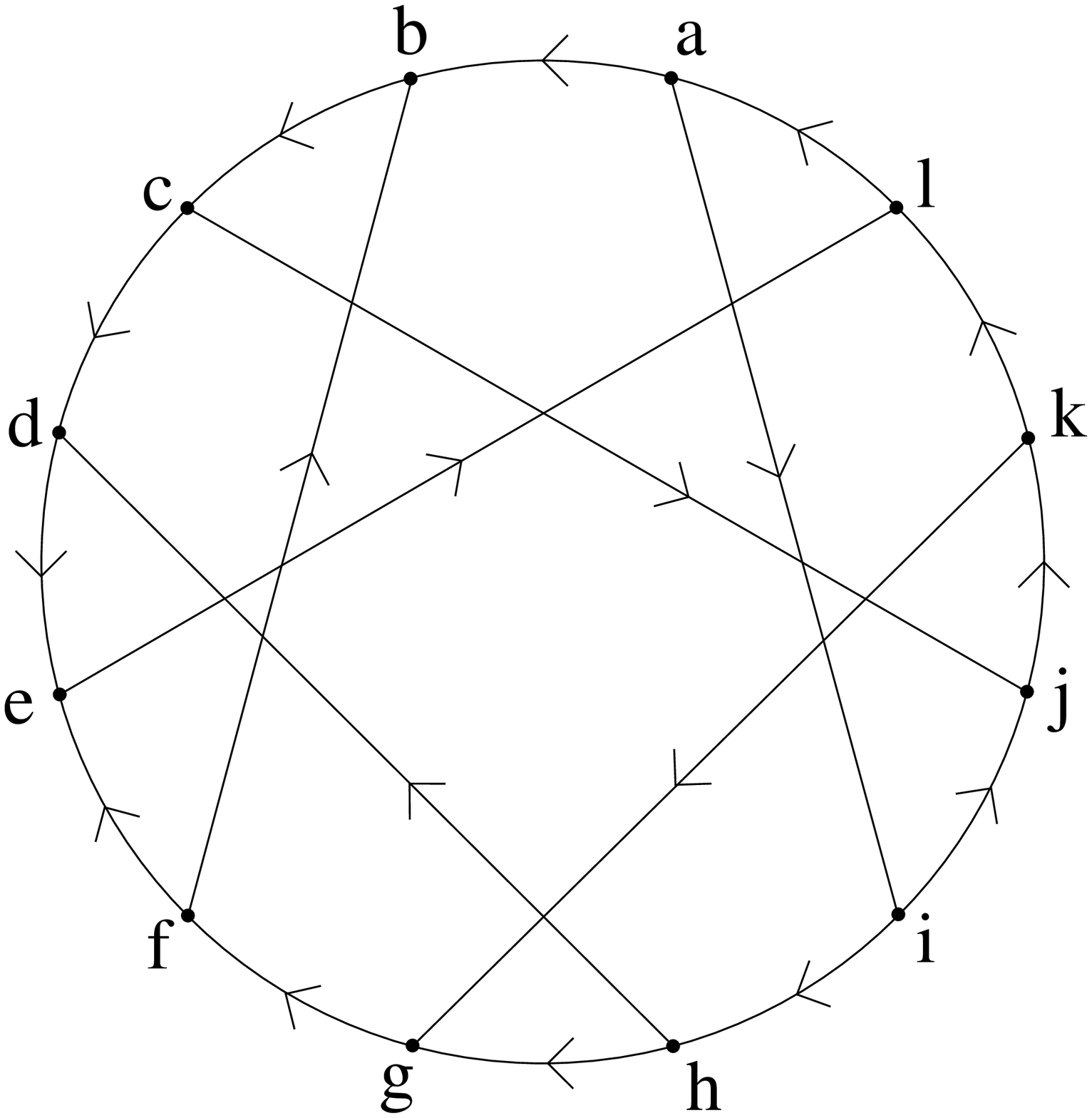}}}%
\hspace{1.5cm} \parbox[b]{1cm}{
\begin{tabular}[b]{cc}
$(a,b)(c,d)(f,e)(h,g)(i,j)(k,l)$&$ +$\\
$(b,c)(d,e)(g,f)(i,h)(j,k)(l,a)$&$ -$\\
$(e,l)(a,b)(c,d)(g,f)(i,h)(j,k)$&$ +$\\
$(h,d)(k,g)(b,c)(f,e)(i,j)(l,a)$&$ -$\\
$(a,i)(f,b)(c,j)(d,e)(h,g)(k,l)$&$ +$\\
$(a,i)(f,b)(c,j)(h,d)(e,l)(k,g)$&$ -$
\end{tabular}
\vspace{1.5cm}}
\end{flushleft}
\caption{The graph $\Gamma_2$.}
\label{gamma2}
\end{figure}

In the statement of our main theorem below, $\Gamma_1$ and $\Gamma_2$ refer to the graphs drawn in Figures~\ref{gamma1} and \ref{gamma2} respectively,
where the arrows are to be ignored. Both graphs are near bipartite, since $\Gamma_1 - \{(f,l),(i,c)\}$ and $\Gamma_2 - \{(f,e),(i,j)\}$
are 1-extendible and bipartite. 
Note that $\Gamma_2$ may be obtained from the Petersen graph by subdividing two fixed edges at a maximum
distance apart and then joining the vertices of degree 2 by an edge. These graphs, like $K_{3,3}$, 
can easily be shown to be non-Pfaffian. Indeed, each graph in Figures~\ref{k33}--\ref{gamma2} is accompanied by a set $S$
of 1-factors such that each edge belongs to just two members of $S$ and $S$ contains an odd number of 1-factors of
each kind of sign under the given orientation. The former property of $S$ implies that the latter is still valid if we change the orientation of a single edge.
Therefore the latter property of $S$ is independent of the orientation and 
consequently the graphs cannot be Pfaffian.

It follows that no even subdivision of these graphs is 
Pfaffian. It is shown in \cite{LiReFi} that a graph $G$ is non-Pfaffian if it has a circuit $X$, of odd length,
such that the graph obtained from $G$ by contracting $X$ to a vertex is non-Pfaffian. 
In general, let us say that a graph $G$ is {\it simply reducible} to a graph $H$ if $G$ has a circuit $X$,
of odd length, such that $H$ is obtained from $G$ by contracting $X$. More generally, we say that $G$ is {\it reducible} to a graph $H$
if for some positive integer $k$ there exist graphs $G_0, G_1, \dots, G_k$ such that $G_0 = G$, $G_k = H$ and, for all $i > 0$, $G_{i-1}$ is simply reducible to $G_i$. 
Thus any graph that is reducible to an even 
subdivision of $K_{3,3}$, $\Gamma_1$ or $\Gamma_2$ is non-Pfaffian. In fact, a graph $G$ must be non-Pfaffian
if it has a subgraph $J$ that is reducible to an even subdivision of $K_{3,3}$, $\Gamma_1$ or $\Gamma_2$
and has the property that $G - VJ$ has a 1-factor. The purpose of this paper is to show 
that the converse of this statement holds for near bipartite graphs. 

\smallskip

\medskip

\medskip

\begin{theo}
\label{main}
A near bipartite graph $G$ is non-Pfaffian if and only if $G$ contains a subgraph $J$, reducible to an even subdivision of
$K_{3,3}$, $\Gamma_1$ or $\Gamma_2$, such that $G - VJ$ has a 1-factor.  
\end{theo}

{\bf Definitions and Notation}. The following definitions and notation are fundamental for this paper. Circuits, non-empty paths and, more generally,
subgraphs with no isolated vertices are determined by
their edge sets, and are therefore identified with them in this paper.
If $X$ is a path or circuit in a graph $G$, then we denote by $VX$ the set of vertices of $X$. If $P$ is a path and $u,v \in VP$,
then we denote by $P[u,v]$ the subpath of $P$ joining $u$ to $v$. If $P[u,v]$ is directed from
$u$ to $v$, then we also write $P(u,v) = P[u,v]$.
If $C$ is a circuit which includes a unique directed path from vertex $u$ to vertex $v$,
then that path is denoted by $C(u,v)$. From time to time we may perform a reorientation of $C$,
that is to say we change the orientation of every oriented edge in $C$.
This directed path is then denoted by $C(v,u)$, or by $P(v,u)$ if it is included in another path $P$. 

A circuit is {\it alternating} with respect to each of two given
1-factors if it is included in their symmetric difference. A circuit that is alternating with respect to a 1-factor $f$
is also said to be {\it $f$-alternating}, or {\it consanguineous} (with respect to $f$). Note that a graph with more than one edge is 1-extendible 
if and only if every edge has an alternating circuit containing it. A path $P$ is {\it alternating}
if every internal vertex of $P$ is incident with an edge of $P \cap f$. An {\it ear} is a path of odd cardinality.

Let $A$ and $B$ be sets of edges in a graph $G$.
Then an {\it $AB$-arc} is a non-empty maximal subpath of $A \cap B$, and an 
{\it $A\overline{B}$-arc} (or a {\it $\overline{B}A$-arc}) is a non-empty maximal subpath of $A - B$. A $G \overline{B}$-arc is also called a $\overline{B}$-arc.  

{\bf 2-Ear Theorem}. Next, let $A$ be an alternating circuit in $G$ and let $H$ be a subgraph of $G$. If there are $n$ $A\overline{H}$-arcs, and each is an ear, then we 
say that $G[EH \cup A]$ is obtained from $H$ by an {\it $n$-ear adjunction}. An {\it ear decomposition} of a 1-extendible graph $G$ is a sequence
$G_0, G_1, \ldots, G_t$ of 1-extendible subgraphs of $G$ such that $G_0$ is isomorphic to $K_2$, $G_t = G$ and,
for each $i > 0$, $G_i$ is obtained from $G_{i-1}$ by a 1-ear or 2-ear adjunction. A theorem of Lov{\'a}sz and Plummer
\cite[Theorem 5.4.6]{LoPl} asserts that every 1-extendible graph has an ear decomposition. It can be stated as follows.

\begin{theo}
\label{2-ear}
Let $f$ be a 1-factor in a 1-extendible graph $G$. Let $H$ be a 1-extendible proper subgraph of $G$ such that 
$EH \neq \emptyset$ and $f \cap EH$ is a 1-factor of $H$. Then $G$ contains an $f$-alternating circuit $A$ that 
admits just one or two $A\overline{H}$-arcs.
\end{theo}

In fact if $G$ is bipartite then it can be shown that only 1-ear adjunctions are necessary.

The idea behind the proof of Theorem~\ref{Bipartite} runs as follows. Clearly we may assume that $G$ is 1-extendible.
Suppose that $G$ is non-Pfaffian. We construct an ear decomposition $G_0, G_1, \dots, G_t$ of $G$. Since $G$ is bipartite, we may assume that,
for each $i > 0$, $G_i$ is obtained from $G_{i-1}$ by the adjunction of a single ear. As $G_0$ is Pfaffian but $G$ is not, there exists a smallest
positive integer $j$ such that $G_j$ is non-Pfaffian. The graph $G_j$ is studied in detail and eventually shown to 
contain $J$.

Theorem~\ref{2-ear} provides a possible way to generalise this argument. If we drop the assumption 
that $G$ is bipartite then, for each $i$, $G_i$ is obtained from $G_{i-1}$ by the adjunction of one or two ears.
In this paper we consider the case where $G_{j-1}$ is bipartite and $G_j$ is obtained from $G_{j-1}$ 
by a $2$-ear adjunction. 

{\bf Idea behind the proof of Theorem~\ref{main}}. We use alternating circuits in preference to 1-factors.
Kasteleyn \cite{Ka} has shown that the 1-factors of a directed graph all have equal sign if and only if all the 
alternating circuits are clockwise odd. (The clockwise parity of a circuit of even length is the parity of the
number of its edges that are directed in agreement with a specified sense.)  Let $G$ be a near bipartite graph which is
minimal with respect to the property of being non-Pfaffian.  Let $e_1$ and $e_2$ be edges of $G$ such that $G - \{e_1, e_2\}$
is bipartite and 1-extendible. By minimality $G - \{e_1, e_2\}$ has a Pfaffian
orientation. Extend this orientation to an orientation of $G$ by orienting $e_1$ and $e_2$ arbitrarily. Since $G$ is
non-Pfaffian, there exist two alternating circuits $A$ and $B$ of opposite clockwise parity. In Theorem~\ref{haupt} we construct alternating circuits 
in $G - \{e_1, e_2\}$ whose union is $EG - \{e_1,e_2\}$ and whose sum (symmetric difference) is $A + B$. This construction is used to 
generate all the non-Pfaffian near bipartite graphs.
The list of non-Pfaffian near bipartite graphs so constructed is infinite. In Sections~3 and 4 we are then able to reduce this list to a finite list by invoking 
the minimality of $G$. In Section~5 we finally show that every graph in this list can be obtained from 
$K_{3,3}$, $\Gamma_1$ or $\Gamma_2$ by means of the operations of reduction and even subdivision. In Section~6 we demonstrate that neither $\Gamma_1$ nor $\Gamma_2$
is reducible to an even subdivision of $K_{3,3}$.            

\section{A structure theorem of minimal non-Pfaffian near bipartite graphs}

In this section we establish that a minimal non-Pfaffian near bipartite graph 
is the union of two alternating circuits $A$ and $B$ 
and two additional paths $S$ and $T$. 
Let $G$ be a near bipartite graph. We may assume that $G$ is minimal with respect to the property of being non-Pfaffian.
To see this point, suppose that $G$ has an edge $e$ such that $G - \{e\}$ is non-Pfaffian and has a subgraph $J$,
reducible to an even subdivision of $K_{3,3}, \Gamma_1$ or $\Gamma_2$, such that $(G - \{e\}) - VJ$ has a 1-factor $f$. Then $f$ is also a 
1-factor of $G - VJ$, and so Theorem~\ref{main} holds also for $G$.

A set $S$ of alternating circuits in a directed graph $H$ is called {\it intractable} if the sum of the circuits in $S$ is empty and an odd number of the members of $S$ are
clockwise even.
The former property implies that the latter is independent of the orientation of $H$. (See Lemma~\ref{orient}.) The following lemma is proved in \cite{Li2}. 
\begin{lem}
A graph is Pfaffian if and only if it has no intractable set of alternating circuits.
\end{lem}

From this result we show that we can assume there to be no vertices of degree 2 in $G$.

\begin{lem}
\label{degree2}
Let $v$ be a vertex of degree $2$ in $G$, and let $G'$ be the graph obtained from $G$ by contracting the edges incident on $v$.
Then $G$ is Pfaffian if and only if $G'$ is Pfaffian.
\end{lem}

{\it Proof}: Let $a$ and $b$ be the edges of $G$ incident on $v$, and let $u$ and $w$ be the vertices adjacent to $v$. 

Suppose there is an intractable set $S$ of alternating circuits in $G$. Then the intersections of the circuits in $S$ with $EG - \{a, b\}$ yield 
an intractable set in $G'$. Conversely, let $S'$ be an intractable set of alternating circuits in $G'$. Let $v'$ be the vertex in $G'$ obtained by 
identifying $u$ and $w$ in $G$. Choose $C' \in S'$. If $v' \notin VC'$, or the edges of $C'$ incident on $v'$ in $G'$ are both 
incident on $u$ in $G$ or both incident on $w$ in $G$, then let $C = C'$; otherwise let $C = C' \cup \{a,b\}$. The set $S$ of such circuits $C$ forms an
intractable set in $G$. (Note that the sum of the circuits in $S$ is a subset of $\{a,b\}$ and therefore empty as it must be a cycle.) \qed 

\smallskip

Let $G$ be a graph with a vertex of degree 2 and let $G'$ be the graph obtained from $G$ by contracting the edges incident on it. Suppose that in $G'$
there is a subgraph $J$, reducible to an even subdivision of $K_{3,3}, \Gamma_1$ or $\Gamma_2$, such that $G' - VJ$ has a 1-factor. Then the same is true for
$G$, for $K_{3,3}$, $\Gamma_1$ and $\Gamma_2$ are cubic and so the converse of the reduction in the lemma gives an even subdivision of each of those graphs.
Therefore we can assume that $G$ has no vertex of degree 2.    

Since $G$ is near bipartite, it is 1-extendible. Moreover there exist edges $e_1$ and $e_2$ 
such that $G - \{e_1, e_2\}$ is bipartite and 1-extendible. We call this graph $H$, and fix a 1-factor $f$ of $H$. Note that $G - \{e_1\}$ is non-bipartite,
for otherwise, since $G$ is non-bipartite, every circuit containing $e_1$ would be of odd length, in contradiction to the fact that $G$ has an alternating 
circuit containing $e_1$.  Similarly $G - \{e_2\}$ is non-bipartite. Consequently any
alternating circuit containing one of $e_1$ and $e_2$ must also contain the other. 

Note that $H$ is Pfaffian, by the minimality of $G$.
Extend a Pfaffian orientation of $H$ to an orientation of $G$ by orienting $e_1$ and $e_2$ arbitrarily. 
We shall henceforth refer to this orientation as our extended Pfaffian orientation of $G$. As $G$ is non-Pfaffian,
it contains a clockwise even alternating circuit $A$. This circuit must contain $e_1$ and $e_2$. There must also be a clockwise odd
alternating circuit $B$ containing $e_1$ and $e_2$, for otherwise a Pfaffian orientation for $G$ could be constructed by
reorienting $e_1$ or $e_2$. The following lemma, which is proved in \cite{LiReFi}, gives information about how $A$ and $B$
can be chosen.

\begin{lem}
\label{AandB}
Let $f$ be a 1-factor in a 1-extendible directed graph $G$. Let $A$ and $B$ be $f$-alternating circuits in $G$, 
of opposite clockwise parity, containing distinct independent edges $e_1$ and $e_2$ such that $e_1 \notin f$
and $e_2 \notin f$. Suppose that $G - \{e_1\}$ and $G - \{e_2\}$ are not bipartite but that $G - \{e_1, e_2\}$ is. 
Then $A \cup B$ includes alternating circuits $X$ and $Y$, of opposite clockwise parity and consanguineous
with respect to some 1-factor that contains neither $e_1$ nor $e_2$, such that there are just one or two $XY$-arcs, each $XY$-arc contains $e_1$ or $e_2$ and 
their union contains both.
\end{lem}

Thus $A$ and $B$ can be chosen so that there are at most two $AB$-arcs. In \cite{LiReFi} the case where 
there is a unique $AB$-arc has been dealt with. We obtained the following theorem.

\begin{theo}
Let $G$ be a 1-extendible graph with 1-factor $f$. Let $e_1$ and $e_2$ be distinct independent edges of $EG - f$
such that neither $G - \{e_1\}$ nor $G - \{e_2\}$ is bipartite but $G - \{e_1, e_2\}$ is bipartite,
Pfaffian and 1-extendible. Suppose there exist $f$-alternating circuits $A$ and $B$, both containing $e_1$ and $e_2$,
such that there is a unique $AB$-arc and $A$ and $B$ have opposite clockwise parity under a Pfaffian orientation of $G - \{e_1, e_2\}$.
Then $G$ has a subgraph $J$, reducible to an even subdivision of $K_{3,3}$, such that $G - VJ$ has a 1-factor. 
\end{theo}

In the present paper, we deal with the remaining case, where for every choice of $A$ and $B$ there are at least
two $AB$-arcs. Henceforth we assume that $A$ and $B$ are chosen so that there are exactly two $AB$-arcs, and therefore exactly two $\overline{A}B$-arcs
and exactly two $A \overline{B}$-arcs. By Lemma \ref{AandB} it may be assumed that one of the $AB$-arcs contains $e_1$ and the other $e_2$.  Let the former arc join
vertices $x_1$ and $x_2$ and the latter vertices $y_1$ and $y_2$. Let $e_1$ join vertices $u_1$ and $u_2$ and let $e_2$ join vertices $v_1$ and $v_2$.
Define $A^* = A - \{e_1\}$, and adjust the notation so that the vertices $u_1$, $x_1$, $y_1$, $v_1$, $v_2$, $y_2$, $x_2$, $u_2$ appear in that order as $A^*$
is traversed from $u_1$ to $u_2$.  

\begin{lem}
One of the $\overline{A} B$-arcs joins $x_2$ to $y_1$ and the other $x_1$ to $y_2$.
\end{lem}

{\it Proof}: Suppose the contrary.
Note that the edges of $f$ incident on $x_1$, $x_2$, $y_1$ and $y_2$, respectively, belong to
$A^*[x_1, u_1] \cup A^*[u_2, x_2] \cup A^*[y_2, y_1]$, since $e_1$ and $e_2$ belong to $AB$-arcs. If an
$\overline{A}B$-arc $X$ were to join $y_1$ to $y_2$ then we should have the contradiction that the circuit $A^*[y_1,y_2] \cup X$ would be of even length
yet contain $e_2$ but not $e_1$. On the other hand, suppose that an $\overline{A}B$-arc $Y$ were to join $x_1$ to $y_1$. 
Let 
$$C = Y \cup A^*[x_1, u_1] \cup \{e_1\} \cup A^*[u_2, y_1].$$
This circuit, an $f$-alternating circuit containing $e_1$ and $e_2$, would have opposite clockwise parity from either $A$ or $B$.
Since there would be a unique $AC$-arc and a unique $BC$-arc, we should have a contradiction to the assumption that there is no choice
for $A$ and $B$ that gives a unique $AB$-arc. \qed

\smallskip

\begin{figure}
\begin{center}
\setlength{\unitlength}{1cm}
\hspace{0.5cm}\mbox{\scalebox{0.35}{%
 \includegraphics{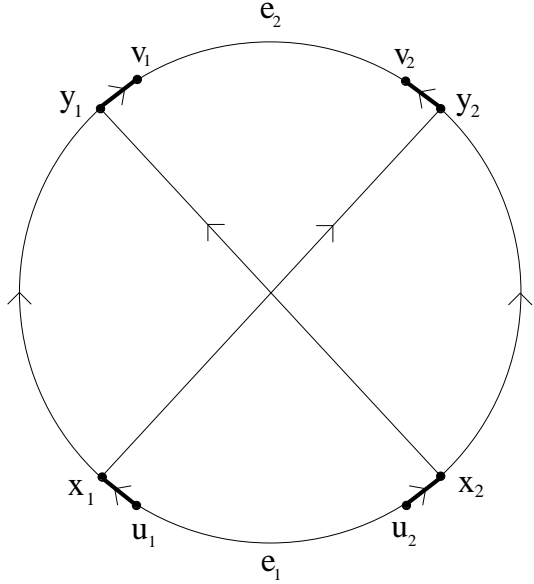}}}%
\end{center}
\caption{A homeomorph of $G[A \cup B]$.}
\label{AB}
\end{figure}

The graph $G[A \cup B]$ is an even subdivision of that shown in Figure~\ref{AB}.
The edges of $f$ are thickened in this and subsequent figures, and in all subsequent figures the graph in question is an even
subdivision of the one portrayed.

For a bipartite graph $K$ with bipartition $\{M, N\}$ and 1-factor $f$ there exists an orientation in which the directed paths and directed circuits
are precisely the $f$-alternating paths and $f$-alternating circuits respectively: orient the edges of $f$ from $M$ to $N$ and the remaining edges from $N$ to $M$.
Then every vertex has indegree 1 or outdegree 1, and every edge joins a vertex of indegree 1 to a vertex of outdegree 1.
It follows that directed paths with an internal vertex in common meet in an edge incident on the vertex.
We call this orientation the {\it reference} orientation for $K$ with respect to $(M,N,f)$.  
Fix such an orientation for $H$ so that $A^*[u_1, v_1]$ is directed from $u_1$ to $v_1$. We refer to this orientation
as our reference orientation for $H$. It follows that
$B \cap EH$ includes a directed path from $u_1$ to $v_2$ and another from $u_2$ to $v_1$, and that 
$A^*[u_2, v_2]$ is directed from $u_2$ to $v_2$. 
The orientation is also indicated in Figure~\ref{AB}.
Henceforth the orientation associated with $H$ will be our reference 
orientation unless an indication to the contrary is given.

Let $f'$ be another 1-factor of $K$. It is shown in \cite{Li1} that the reference orientation for $K$ with respect to $(M,N,f')$ is obtained from the reference orientation with
respect to $(M,N,f)$ by reorienting the circuits included in $f + f'$. This fact is used implicitly later on to justify reorientations
of $f$-alternating circuits.

The following lemma is a standard result. (See \cite{Li1}.)

\begin{lem}
\label{standard}
Let $G$ be a directed graph such that each edge has a directed circuit containing it. Then for every 
$a,b \in VG$, there exists a directed path from $a$ to $b$. 
\end{lem}

\begin{figure}
\begin{center}
\setlength{\unitlength}{1cm}
\hspace{0.5cm}\mbox{\scalebox{0.35}{%
 \includegraphics{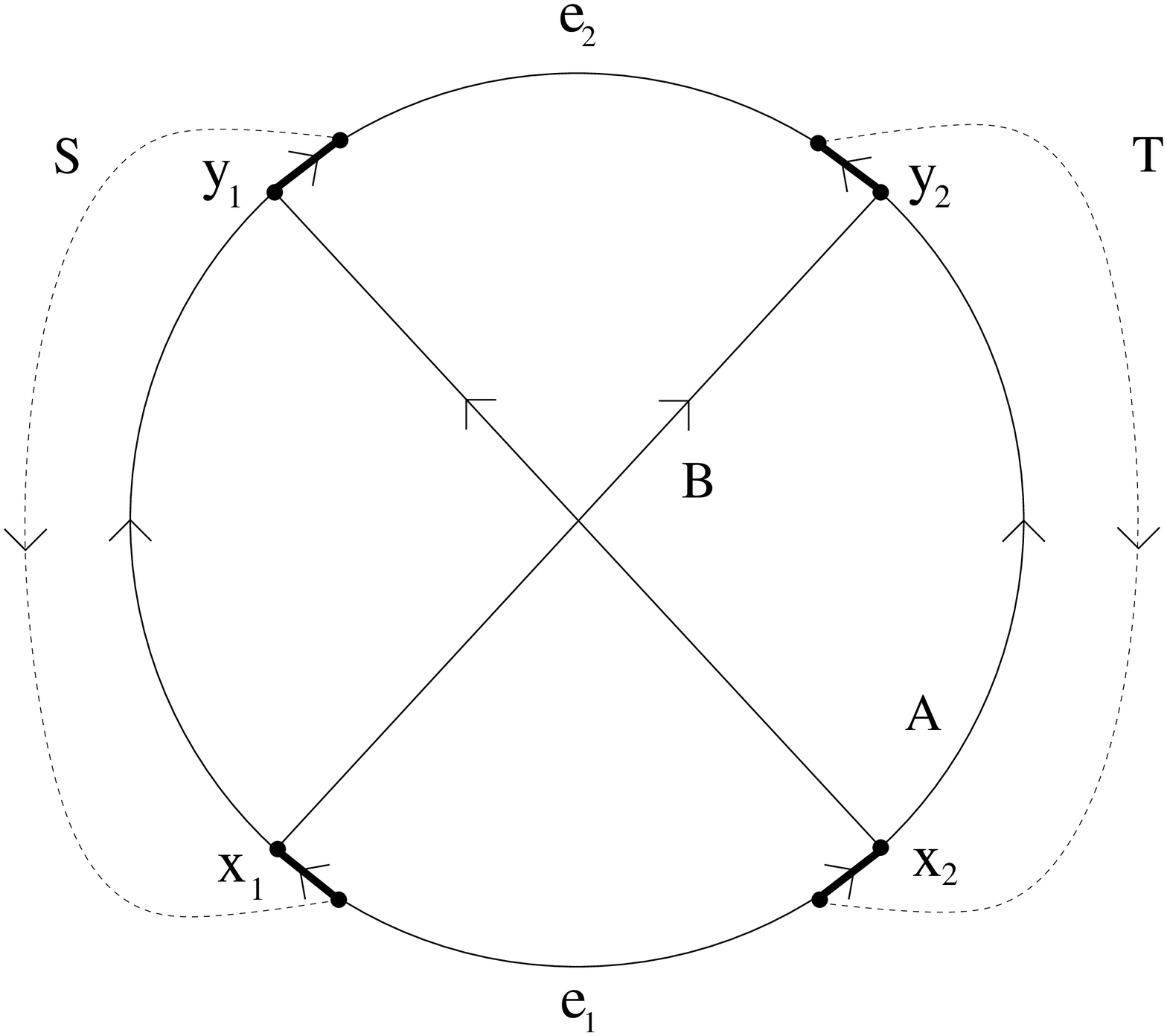}}}%
\end{center}
\caption{$A \cup B \cup S \cup T$}
\label{hauptfig}
\end{figure}

We may apply this lemma to $H$, since every edge of the 1-extendible graph $H$ must belong to a directed
circuit. Thus there is a directed path $S$ from $y_1$ to $x_1$ and a directed path $T$ from $y_2$ to $x_2$. (See Figure~\ref{hauptfig};
a dotted line in this and subsequent figures stands for a directed path, which can have intersections with the rest of the graph
that are not indicated.)
We now aim to show that $EG = A \cup B \cup S \cup T$ in Theorem~\ref{haupt}. 
To this end we introduce the following three lemmas.

\begin{lem}
\label{seq}
Let $(a_1, a_2, \dots, a_n)$ be a sequence of edges in a directed graph $G$ in which each vertex has indegree $1$ or outdegree $1$. Suppose that for all $i > 1$ the
origin of $a_i$ is the terminus of $a_{i-1}$ and that the origin of $a_1$ has outdegree $1$ and the terminus of $a_n$ has indgree $1$.
Then there exist a directed path $P$ from the origin of $a_1$ to the terminus of $a_n$ and directed circuits 
$C_1, C_2, \dots, C_k$ such that 
  $$\sum_{i=1}^n \{a_i\} = P + \sum_{i=1}^k C_i$$
and 
$$ P \cup \bigcup_{i=1}^k C_i = \{a_1,a_2,\dots,a_n\}.$$ 
\end{lem}

{\it Proof}:  
We use induction on the number $r$ of repetitions of edges in the sequence $L = (a_1, a_2, \dots, a_n)$.
If $r = 0$ then $P = \{a_1, a_2, \dots, a_n\}$ and $k = 0$, since each vertex has indegree 1 
or outdegree 1, the origin of $a_1$ has outdegree $1$ and 
the terminus of $a_n$ has indegree $1$.  
Now suppose that $r > 0$ and that the lemma holds whenever
the number of repetitions of edges is less than $r$. Let $a$ be the edge in $L$ that is repeated first.
The part of $L$ between the first two occurrences of $a$ has no edges repeated within it.
Therefore $a$ and the edges between the first two occurrences of $a$ form a directed circuit $C$. We now modify $L$ 
by removing all the edges from the first occurrence of $a$ to the edge immediately before the second occurrence of $a$. This 
modified sequence $L' = (a_1', a_2', \dots, a_m')$ has fewer repetitions of edges. 
Moreover the origin of $a_1'$ is that of $a_1$, the terminus of $a_m'$ is that of $a_n$, and for all $i > 1$
the origin of $a_i'$ is the terminus of $a_{i-1}'$. Therefore the inductive hypothesis may be applied to $L'$, and the
result follows from the equation 
  $$ \sum_{i=1}^n \{a_i\} = \sum_{i=1}^m \{a_i'\} + C. $$ \qed

\begin{lem}
\label{circuits}
Let $H$ be a directed graph, let $P$ be a directed path from vertex $x$ to vertex $y$ and let $Q$
be a directed path from $y$ to $x$. Then there are directed circuits $C_1, C_2, \dots, C_k$ such that 
   $$\sum_{i=1}^k C_i = P + Q$$
and 
$$ \bigcup_{i=1}^k C_i = P \cup Q. $$ 
\end{lem}

{\it Proof}: We use induction on $n = |VP \cap VQ| \ge 2$. If $n = 2$, then $\{P \cup Q\}$ is the required
set of directed circuits.

\begin{figure}
\begin{center}
\setlength{\unitlength}{1cm}
\hspace{0.5cm}\mbox{\scalebox{0.35}{%
 \includegraphics{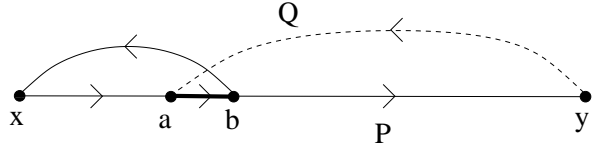}}}%
\end{center}
\caption{P and Q in Lemma~\ref{circuits}.}
\label{lincirc}
\end{figure}

Let $n > 2$ and suppose the lemma holds whenever $|VP \cap VQ| < n$. Let $b$ be the last vertex of 
$VQ - \{x\}$ that is in $VP$, and let $a$ be the last vertex of $VQ$ incident with an edge of $Q(y,b) - P$. (See Figure~\ref{lincirc}.) By the inductive
hypothesis there exist circuits $C_1, C_2, \dots, C_{k-1}$ such that 
  $$ \sum_{i=1}^{k-1} C_i = P(a,y) + Q(y,a).$$
The required set of circuits is
  $$ \{C_1, C_2, \dots, C_{k-1}, P(x,b) \cup Q(b,x)\},$$
since
  \begin{eqnarray*}
   && P(a,y) + Q(y,a) + (P(x,b) \cup Q(b,x)) \\ 
   &=& P(a,y) + Q(y,a) + P(x,b) + Q(b,x) \\
   &=&  P(a,y) + Q(y,a) + P(x,a) + P(a,b) + Q(b,x) \\ 
   &=& P(x,a) + P(a,y) + Q(y,a) + Q(a,b) + Q(b,x) \\
   &=& P + Q  
  \end{eqnarray*}
as $P(a,b) = Q(a,b)$. \qed

\begin{lem}
\label{orient}
Let ${\mathcal C}$ be a set of circuits of even length and empty sum in a directed graph $G$. Then the parity of the number
of clockwise even members of ${\mathcal C}$ is independent of the orientation of $G$.  
\end{lem}

{\it Proof}: A change of orientation can be effected by changing orientations of edges one at a time. Each such change leaves
unaltered the parity of the number of clockwise even circuits in ${\mathcal C}$. \qed

\begin{theo}
\label{haupt}
Let $G$ be a minimal non-Pfaffian near bipartite graph. Let $e_1$ and $e_2$ be edges such that $G - \{e_1, e_2\}$ is bipartite and 1-extendible.
Let $H = G - \{e_1, e_2\}$, and let $f$ be a 1-factor of $H$. Let $A$ and $B$ be $f$-alternating circuits in $G$ of opposite clockwise parity.
Suppose that there are exactly two $AB$-arcs, one containing $e_1$ and the other $e_2$. Let the former arc join vertices $x_1$ and $x_2$ and the
latter vertices $y_1$ and $y_2$. Let $A^* = A - \{e_1\}$, and suppose the vertices $x_1$, $y_1$, $y_2$, $x_2$
appear in that order when $A^*$ is traced from $x_1$ to $x_2$. Let $H$ be given its reference orientation with respect to $f$ such that $A^*[x_1, y_1]$
is a directed path from $x_1$ to $y_1$. Let $S'$
be a directed path from $y_i$ to $x_j$ and $T'$ a directed path from $y_{3-i}$ to $x_{3-j}$, where $i, j \in \{1,2\}$. Then $G$ is induced by
$A \cup B \cup S' \cup T'$.    
\end{theo}

{\it Proof}: Without loss of generality we take $i=j=1$. We first show, by using Lemmas~\ref{seq} and \ref{circuits}, that we can write $A + B$ as a sum of directed circuits
included in $A \cup B \cup S' \cup T'$.
In order to verify this claim, first we apply Lemma~\ref{circuits} to the directed paths $A(x_1,y_1)$ and $S'$. Let
  \begin{equation} \label{eq1} A(x_1,y_1) + S' = \sum_{V \in {\mathcal V}} V, \end{equation}
where ${\mathcal V}$ is a set of directed circuits included in $A(x_1,y_1) \cup S'$. Similarly let 
  \begin{equation} \label{eq2} A(x_2,y_2) + T' = \sum_{W \in {\mathcal W}} W, \end{equation}
where ${\mathcal W}$ is a set of directed circuits included in $A(x_2,y_2) \cup T'$. 
Now consider the sequence $L$ of edges formed by the edges in the directed path $S'$ followed by those in the directed path $B(x_1,y_2)$.
We apply Lemma~\ref{seq} to $L$ to write 
  \begin{equation}  \label{eq3} S' + B(x_1,y_2) = P + \sum_{X \in {\mathcal X}} X, \end{equation}
where $P$ is a directed path from $y_1$ to $y_2$ included in $S' \cup B(x_1,y_2)$ and ${\mathcal X}$ is a set of directed circuits
included in $S' \cup B(x_1,y_2)$. Similarly we have
  \begin{equation}  \label{eq4} T' + B(x_2,y_1) = Q + \sum_{Y \in {\mathcal Y}} Y, \end{equation}
where $Q$ is a directed path from $y_2$ to $y_1$ included in $T' \cup B(x_2,y_1)$ and ${\mathcal Y}$ is a set of directed circuits included in $T' \cup B(x_2,y_1)$.
Now we apply Lemma~\ref{circuits} to $P$ and $Q$ to obtain
  \begin{equation}  \label{eq5} P + Q = \sum_{Z \in {\mathcal Z}} Z, \end{equation}
where ${\mathcal Z}$ is a set of directed circuits included in $P \cup Q$. Let ${\mathcal C} = {\mathcal V} + {\mathcal W} + {\mathcal X} + {\mathcal Y} + {\mathcal Z}.$
Adding equations~\eqref{eq1}--\eqref{eq5} we obtain
\begin{eqnarray*}
&& A(x_1,y_1) + A(x_2,y_2) + B(x_1,y_2) + B(x_2,y_1) \\
& = & \sum_{V \in {\mathcal V}} V + \sum_{W \in {\mathcal W}} W + \sum_{X \in {\mathcal X}} X + \sum_{Y \in {\mathcal Y}} Y + \sum_{Z \in {\mathcal Z}} Z \\
& = & \sum_{C \in {\mathcal C}} C.  
\end{eqnarray*}
Since the left hand side is $A + B$, ${\mathcal C}$ is the required set of circuits.

In our extended Pfaffian orientation of $G$, $A$ is the only clockwise even circuit in ${\mathcal C} \cup \{A,B\}$.
Therefore by Lemma~\ref{orient} an odd number of circuits in ${\mathcal C} \cup \{A,B\}$ are clockwise even for any orientation of $G$.
But if there were a Pfaffian orientation of $G[A \cup B \cup \bigcup_{C \in {\mathcal C}} C]$ then every circuit in ${\mathcal C} \cup \{A,B\}$ would be clockwise odd because
they are $f$-alternating. Therefore the graph $G[A \cup B \cup \bigcup_{C \in {\mathcal C}} C]$
is non-Pfaffian, and so $G[A \cup B \cup S' \cup T']$ is non-Pfaffian. By the minimality of $G$, we deduce
that
 $$G[A \cup B \cup S' \cup T'] = G.$$  \qed

\smallskip

Applying this theorem to $S$ and $T$ we find that $G=G[A \cup B \cup S \cup T]$. In fact we chose $S$ and $T$ to satisfy the 
following definition.

\begin{defi}
Let $i,j \in \{1,2\}$ and let $P$ be a directed path from $y_i$ to $x_j$ under our reference orientation. 
We say that $P$ is {\it minimal} if for every edge
$e \in P - (A \cup B)$ there is no directed path from $y_i$ to $x_j$ included in 
$(A \cup B \cup P) - \{e\}$. 
\end{defi}

It is clear
that $S$ and $T$ may be assumed to be minimal. 

Let $P$ be a directed path from vertex $x$ to vertex $y$. Let $P_1$ and $P_2$ be disjoint subpaths of $P$ such that each edge of $P_1$
is closer to $x$ in $P$ than is any edge of $P_2$. In this case we write $P_1 <_P P_2$. If $P_1 = \{a_1\}$ and $P_2 = \{a_2\}$, then we write $a_1 <_P a_2$ instead
of $\{a_1\} <_P \{a_2\}$. A similar notation is used for vertices in $P$. 

The next lemma is obvious.
\begin{lem}
\label{unique}
Let $u$ and $v$ be vertices in $G[A \cup B]$. Under our reference orientation for $H$ there exists at most one directed path in $G[A \cup B]$ from $u$ to $v$. 
\end{lem}

This lemma is used in the proof of the following lemma.

\begin{figure}
\begin{center}
\setlength{\unitlength}{1cm}
\mbox{\scalebox{0.35}{%
 \includegraphics{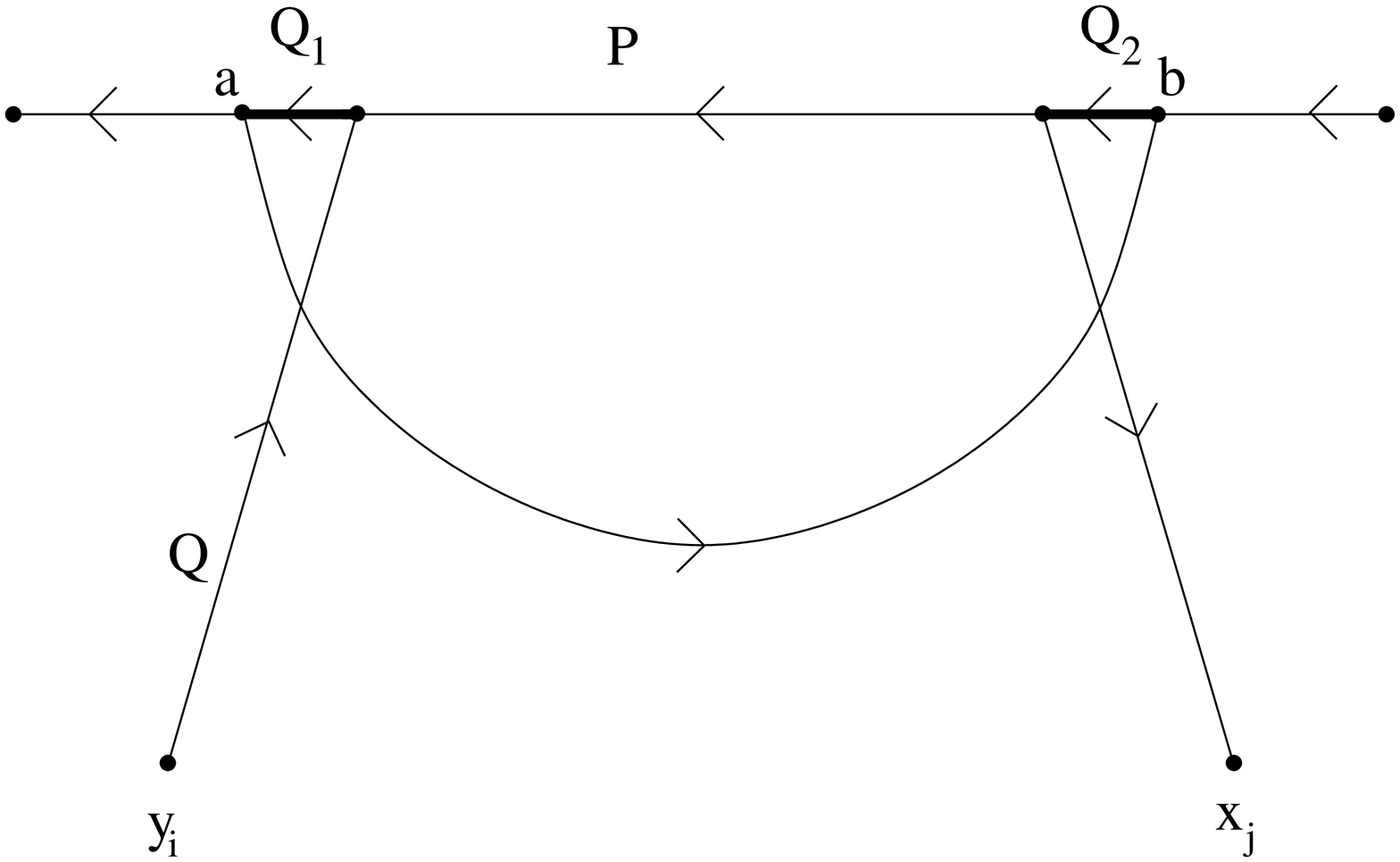}}}%
\end{center}
\caption{A directed minimal path $P$ from $y_i$ to $x_j$.}
\label{minimal1}
\end{figure}

\begin{lem}
\label{minimal}
(a) Let $Q$ be a minimal directed path from $y_i$ to $x_j$ 
and let $P$ be a directed path included in $A \cup B$.
Let $Q_1$ and $Q_2$ be distinct $QP$-arcs. Then $Q_1 <_Q Q_2$ if and only if $Q_2 <_P Q_1$. (See Figure~\ref{minimal1}.)\\
(b) Conversely let $Q'$ be a directed path from $y_i$ to $x_j$. If for 
every directed path $P'$ in $A \cup B$ and every pair of distinct $Q'P'$-arcs $Q'_1$ and $Q'_2$ 
we have  $Q'_1 <_{Q'} Q'_2$ if and only if $Q'_2 <_{P'} Q'_1$, then $Q'$ is minimal.
\end{lem}

{\it Proof}: (a) It suffices to show that if $Q_1 <_Q Q_2$ then $Q_2 <_P Q_1$. 
Assume the contrary, that is $Q_1 <_Q Q_2$ and $Q_1 <_P Q_2$.  
Let $a$ be the terminus of $Q_1$
and $b$ the origin of $Q_2$. By assumption $P(a,b)$ and $Q(a,b)$ exist, 
but it are not equal. By Lemma~\ref{unique} there is an edge 
$e \in Q(a,b) - (A \cup B)$. The set $Q(y_i, a) \cup P(a,b) \cup Q(b, x_j)$ includes a directed path from $y_i$ to $x_j$.
This path is included in $(A \cup B \cup Q) - \{e\}$, in contradiction to the minimality of $Q$. 

(b) Conversely, assume that $Q'$ is not minimal. Choose $e \in Q' - (A \cup B)$ so that there exists a directed path $Q$ from $y_i$ to $x_j$
in $(A \cup B \cup Q') - \{e\}$. Let $u$ be the last vertex in $Q$ that is also in $Q'$ and satisfies $Q(y_i, u) = Q'(y_i, u)$. Let $v$ be the first vertex 
in $VQ(u, x_j) - \{u\}$ that is also in $Q'$. Then $u <_Q v$, $u <_{Q'} v$ 
and $Q(u, v) \cap Q' = \emptyset$. Hence $Q(u,v) \subseteq A \cup B$ and so $Q(u,v)$ is included 
in a maximal directed path $P'$ included in $A \cup B$ such that $u <_{P'} v$. 
Let $Q'_1$ be the $Q'P'$-arc that includes $u$ and $Q'_2$ the $Q'P'$-arc that includes
$v$. Then $Q'_1$ and $Q'_2$ are distinct, $Q'_1 <_{Q'} Q'_2$ but $Q'_1 <_{P'} Q'_2$. \qed

\section{Forbidden $\overline{A \cup B}$-arcs}

In this section we rule out certain directed $\overline{A \cup B}$-arcs.
For that purpose we need the following technical lemma.

\begin{lem}
\label{include}
Let $P$ be a directed $\overline{A \cup B}$-arc. Then there exist $i,j \in \{1,2\}$ such that $P$ is included in a minimal path
$Q$ directed from $y_i$ to $x_j$.
\end{lem}

{\it Proof}: If $S \cap T = \emptyset$, then 
$S$ and $T$ are vertex disjoint and, since $P \subseteq (S \cup T) - (A \cup B)$ by Theorem~\ref{haupt}, it follows that $P \subseteq S$
or $P \subseteq T$. Assume therefore that
$S \cap T \neq \emptyset$. Let $a$ and $b$ be, respectively, the first and last vertices of $T$ that are also in $S$.
It follows from Theorem~\ref{haupt} that
\begin{equation}
\label{G} 
   G = G[A \cup B \cup S \cup T(y_2,a) \cup T(b,x_2)], 
\end{equation}
as there exists a directed path from $y_2$ to $x_2$ included in
  $$ T(y_2,a) \cup S(a,x_1) \cup A(x_1,y_1) \cup S(y_1,b) \cup T(b,x_2). $$

We observe from \eqref{G} that any vertex of degree 3 and not in $VA \cup VB$ must be either $a$ or $b$.
It follows that if either $a$ or $b$ were not an internal vertex of $P$, then $P$ would be included in $S$ or $T$, since the 
edges of $f$ incident on $a$ or $b$ are in $S \cap T$. In this case we could choose
$Q$ to be $S$ or $T$. We therefore assume that $a$ and $b$ are the internal vertices of $P$, since $G$ has no vertices of degree $2$.

Let $u$ and $v$ be, respectively, the origin and terminus of $P$.
If $b <_P a$, then $P(u,b) \cup P(a,v) \subseteq S \cap T$, and so $P \subseteq S$
or $P \subseteq T$ according to whether $P(b,a) \subseteq S$ or $P(b,a) \subseteq T$. Therefore we can assume that $a <_P b$.
Then $S(a,b) = T(a,b)$. Moreover $P(u,a)$ is included in $S$ or $T$, and similarly for $P(b,v)$. Without loss of generality we
assume that $P(u,a) \subseteq S$. If $P(b,v) \subseteq S$, then we take $Q = S$. Suppose therefore that $P(b,v) \subseteq T$.
In this case, we take $Q = S(y_1,b) \cup T(b,x_2)$. 

It remains to show that $Q$ is minimal. Suppose not. Then there exists $e \in Q - (A \cup B)$ such that there is a directed path $Q'$ 
from $y_1$ to $x_2$ included in $(A \cup B \cup Q) - \{e\}$. Define $R = T(y_2, b) \cup S(b, x_1)$. Then by Theorem~\ref{haupt} we have
$G = G[A \cup B \cup R \cup Q']$. If $e \in S(y_1, a) \cup T(b, x_2)$ then we have the contradiction that $e \notin A \cup B \cup R \cup Q'$.
Therefore we suppose that $e \in Q(a,b)$. Then $P \cap Q' = \emptyset$, and we have the contradiction that 
  $$ (S(u,a) \cup T(b,v)) \cap (A \cup B \cup R \cup Q') = \emptyset. $$ \qed

\medskip

The next lemma appeared in \cite{LiReFi}.
\begin{lem}
\label{parity}
Let $A_1, A_2$ be $f$-alternating circuits in a directed graph $G$ with 1-factor $f$. Then $A_1$ and $A_2$ are of opposite 
clockwise parity if and only if $A_1 + A_2$ includes an odd number of clockwise even alternating circuits.
\end{lem}

\begin{cor}
\label{A+B}
The sum $A + B$ is a clockwise even circuit under our extended Pfaffian 
orientation.
\end{cor}

{\it Proof}: Note that
  $$ A + B = A(x_1, y_1) \cup B(x_2, y_1) \cup A(x_2, y_2) \cup B(x_1, y_2), $$
which is a circuit. Since $A$ and $B$ are of opposite clockwise parity, the result follows from Lemma~\ref{parity}. \qed

\smallskip

\begin{lem}
\label{long}
Let $P$ be a directed path included in $A \cup B$ such that no internal vertex of $P$ is in $\{x_1,x_2,y_1,y_2\}$. Then there does not exist a directed
$\overline{A \cup B}$-arc joining vertices in $VP$. 
\end{lem}    

\begin{figure}
\begin{center}
\setlength{\unitlength}{1cm}
\mbox{\scalebox{0.35}{%
 \includegraphics{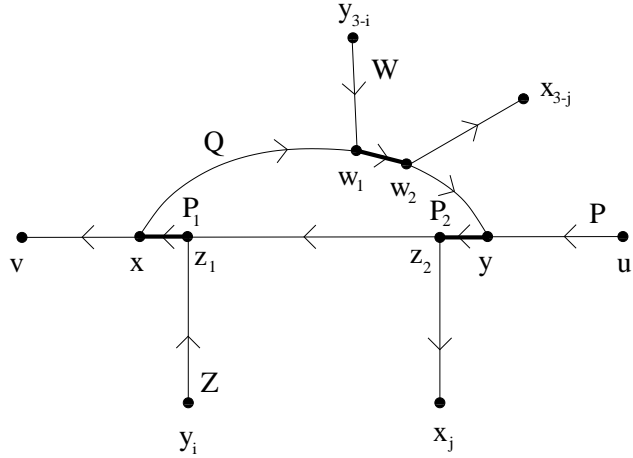}}}%
\end{center}
\caption{The situation in Lemma~\ref{long}.}
\label{longfig}
\end{figure}

{\it Proof}: Suppose there exists a directed $\overline{A \cup B}$-arc $Q$ from $x \in VP$ to $y \in VP$. Then, by Lemma~\ref{include},
for some $i,j \in \{1,2\}$ there exists a directed minimal path $Z$ from $y_i$ to $x_j$ that includes $Q$. By Lemma~\ref{standard} we may choose 
a directed path $W$ from $y_{3-i}$ to $x_{3-j}$; thus $G = G[A \cup B \cup Z \cup W]$ by Theorem~\ref{haupt}.
There exist a $ZP$-arc $P_1$ with terminus $x$ and a $ZP$-arc $P_2 \neq P_1$ with origin $y$. Let $z_1$ be the origin of $P_1$ and $z_2$ the terminus
of $P_2$. (See Figure~\ref{longfig}.)
Since $x <_Z y$ we have $P_1 <_Z P_2$, so that $P_2 <_P P_1$ by Lemma~\ref{minimal}. Therefore $z_2 <_P z_1$. 

Let $C$ be the circuit $Q \cup P(y,x)$. First we show that we may assume there to be at most one $CW$-arc. Suppose there are two such arcs,
$W_1$ and $W_2$, where $W_1 <_W W_2$. Let $a$ be the terminus of $W_1$ and $b$ the origin of $W_2$. Let $W^*$ be a directed path
from $y_{3-i}$ to $x_{3-j}$ included in
   $$W(y_{3-i}, a) \cup C(a, b) \cup W(b, x_{3-j}). $$
The number of $C W^*$-arcs is less than the number of $C W$-arcs. By repeating the argument if necessary and appealing
to the finiteness of $G$, we may therefore assume that there is at most one $CW$-arc. If such an arc exists, let its
origin be $w_1$ and its terminus $w_2$. (See Figure \ref{longfig}.) We also note that there is a unique $CZ$-arc, by the minimality of $Z$.

Let $f'$ be the 1-factor $f + C$, and let $A' = A + C$, $B' = B + C$, $Z' = Z + C$ and $W' = W + C$.
Now we show that
\begin{equation}
\label{A'B'Z'W'} 
       G = G[A' \cup B' \cup Z' \cup W'].  
\end{equation}
By Lemma~\ref{parity}, $A'$ and $B'$ are $f'$-alternating circuits containing $e_1$ and $e_2$ 
of opposite clockwise parity, since $A$ and $B$ are of opposite clockwise parity
and $A' + B' = A + B$. Moreover there are exactly two $A'B'$-arcs and 
the vertices of degree 3 in $G[A' \cup B']$ are the
same as those in $G[A \cup B]$, since 
 $$\{x_1, x_2, y_1, y_2\} \cap (VC - \{x,y\}) = \emptyset.$$ 
In addition $Z'$ would become a directed path from 
$y_i$ to $x_j$ if $C$ were reoriented, and a similar statement holds for $W'$. Thus \eqref{A'B'Z'W'} holds, by Theorem~\ref{haupt} with $f$, $A$ and $B$
replaced by $f'$, $A'$ and $B'$ respectively.

Now we observe that 
  $$ (A' \cup B') \cap P(y,x) = \emptyset $$
since $$\{x_1, x_2, y_1, y_2\} \cap (VP - \{x,y\}) = \emptyset.$$ Note that
  $$ Z' = Z(y_i, z_1) \cup P(z_2, z_1) \cup Z(z_2, x_j). $$
Therefore
  $$ Z' \cap (P_1 \cup P_2) = \emptyset. $$
Thus $ (A' \cup B' \cup Z') \cap (P_1 \cup P_2) = \emptyset$, and so $ P_1 \cup P_2 \subseteq W'$ by~\eqref{A'B'Z'W'}.
We deduce that $w_1$ and $w_2$ exist. 

Next we show that either $w_1, w_2 \in VP(z_2,z_1)$ and $z_2 <_P w_1 <_P w_2 <_P z_1$, or $w_1, w_2 \in VQ$ and $w_1 <_Q w_2$. First,
  $$ W' = W(y_{3-i}, w_1) \cup C(w_2, w_1) \cup W(w_2, x_{3-j}). $$
Hence
  $$ P_1 \cup P_2 \subseteq C(w_2, w_1), $$
and the desired conclusion follows.

{\it Case 1}: Suppose first that $w_1, w_2 \in VP(z_2,z_1)$ and $z_2 <_P w_1 <_P w_2 <_P z_1$. 
After reorientation of $C$, let $X$ be a directed path from $y_i$ to $x_{3-j}$ included in 
  $$ Z(y_i,z_1) \cup C(z_1,w_2) \cup W(w_2,x_{3-j}) $$
and let $Y$ be a directed path from $y_{3-i}$ to $x_j$ included in
  $$W(y_{3-i},w_1) \cup C(w_1, z_2) \cup Z(z_2,x_j). $$
Thus
$G = G[A' \cup B' \cup X \cup Y]$. We now have the contradiction that 
  $$ (P_1 \cup C(w_2, w_1) \cup P_2) \cap (A' \cup B' \cup X \cup Y) = \emptyset. $$

{\it Case 2}: Suppose on the other hand that $w_1, w_2 \in VQ$ and $w_1 <_Q w_2$. Without loss of generality we may assume that $P$ is a maximal
directed path in $A \cup B$ such that no internal vertex is in $\{x_1, x_2, y_1, y_2\}$. Let $P$
be directed from vertex $u$ to vertex $v$. Thus $u \in \{u_1, u_2, x_1, x_2, y_1, y_2\}$ and $v \in \{x_1, x_2, y_1, y_2, v_1, v_2\}$.

First we show that $u \neq y_i$. If $u = y_i$ then we observe that there is a directed path $Z^*$ from $y_i$
to $x_j$ included in $P(u, z_2) \cup Z(z_2, x_j)$. Therefore 
  $$ (A \cup B \cup Z^* \cup W) \cap (C(x, w_1) \cup C(w_2, y)) = \emptyset. $$
By Theorem~\ref{haupt} we have the contradiction that $G = G[A \cup B \cup Z^* \cup W]$. Thus $u \neq y_i$.
A similar argument, with $Z^*$ included in $Z(y_i, z_1) \cup P(z_1, v)$, shows that $v \neq x_j$.

Next we show that $u \neq y_{3-i}$. Otherwise we define $W^*$ to be a directed path from $y_{3-i}$ to $x_{3-j}$ included in 
  $$ P(u, x) \cup C(x, w_2) \cup W(w_2, x_{3-j}). $$
The union $A \cup B \cup Z \cup W^*$ does not contain the edge of $W - Q$ incident on $w_1$. 
This result contradicts Theorem~\ref{haupt}, since
$G = G[A \cup B \cup Z \cup W^*]$. Thus $u \neq y_{3-i}$. A similar argument, with $W^*$ included in 
  $$ W(y_{3-i}, w_1) \cup C(w_1, y) \cup P(y, v), $$
shows that $v \neq x_{3-j}$.

Now we show that $u \neq x_{3-j}$. Otherwise we define $Z^*$ as a directed path from $y_i$ to $x_j$ included in 
  $$ Z(y_i, w_2) \cup W(w_2, u) \cup P(u, z_2) \cup Z(z_2, x_j). $$
Then 
$$(A \cup B \cup Z^* \cup W) \cap C(w_2, y) = \emptyset,$$ 
in contradiction to Theorem~\ref{haupt} since $G = G[A \cup B \cup Z^* \cup W]$.
A similar argument shows that $v \neq y_{3-i}$.

Since $v \notin \{x_1,x_2\}$ we have $u \notin \{u_1,u_2\}$. We conclude that $u = x_j$, and 
similarly $v = y_i$. Define $P' = P + C$. Remember that $W'$ is the only member of $\{A', B', W', Z'\}$ meeting 
$P_1 \cup P_2$. If $C$ is reoriented then there is a directed path $W''$ from $y_{3-i}$ to $x_{3-j}$ included in
  $$ W'(y_{3-i}, w_1) \cup P'(w_1, v) \cup Z'(v, u) \cup P'(u, w_2) \cup W'(w_2, x_{3-j}). $$
Thus 
  $$W'' \cap (P_1 \cup P_2) = \emptyset. $$
We now have a contradiction, since $G = G[A' \cup B' \cup Z' \cup W'']$. \qed 

\begin{lem}
\label{afterlong}
There is no directed $\overline{A \cup B}$-arc joining vertices in distinct $A \overline{B}$-arcs, or in distinct $B \overline{A}$-arcs. 
\end{lem}

\begin{figure}
\begin{center}
\setlength{\unitlength}{1cm}
\mbox{\scalebox{0.35}{%
 \includegraphics{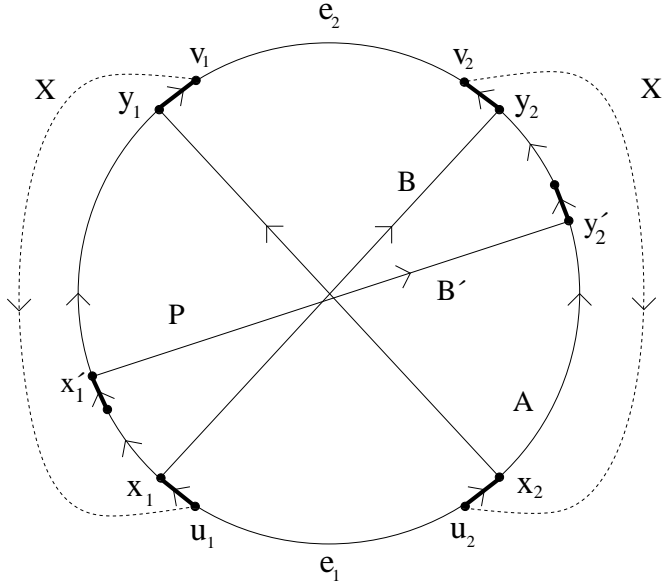}}}%
\end{center}
\caption{$P$ is a  directed arc joining a vertex in $A(x_1,y_1)$ to a vertex in $A(x_2,y_2)$.}
\label{AbBtoAbB}
\end{figure}

{\it Proof}: In view of the symmetry between $A$ and $B$ it suffices to prove that no $\overline{A \cup B}$-arc is directed from a vertex in $A(x_1, y_1)$
to a vertex in $A(x_2, y_2)$. Suppose such an arc $P$ exists, joining a vertex $x_1' \in VA(x_1, y_1)$ to a vertex $y_2' \in VA(x_2, y_2)$ .
(See Figure~\ref{AbBtoAbB}.)
Let 
  $$ B' = A(u_1, x_1') \cup P \cup A(y_2', v_2) \cup \{e_2\} \cup B(u_2, v_1) \cup \{e_1\}. $$
This is an $f$-alternating circuit containing $e_1$ and $e_2$. Observe that $B(x_1, y_2)$ is the only $B \overline{B'}$-arc. It follows
that $B$ and $B'$ have the same clockwise parity, for otherwise $A$ and $B$ could have been chosen to have a unique $AB$-arc.
Henceforth $B'$ will play the r\^ole previously assumed by $B$. The circuit
$A$ will play the same r\^ole as before, but we define $x_2' = x_2$ and $y_1' = y_1$.
Note that there are exactly two $AB'$-arcs, one containing $e_1$ and the other containing 
$e_2$ and that $B(x_1, y_2)$ is an $\overline{A \cup B'}$-arc. Therefore, by Lemma~\ref{include}, for some $i, j \in \{1,2\}$
there exists a directed minimal path $X$ from $y_i'$ to $x_j'$ including $B(x_1, y_2)$. 

We now show that $i=1$ and $j=2$.  
Included in the set $X(y_i', x_1) \cup A(x_1, x_1')$ is a directed path $W$
from $y_i'$ to $x_1'$. This path is included in $(A \cup B \cup X) - B(x_1, y_2)$, in contradiction to the minimality of $X$ if $j=1$.
Therefore $j = 2$.
Similarly, included in the set $A(y_2', y_2) \cup X(y_2, x_2')$ is a directed
path $Z$ from $y_2'$ to $x_2'$. This path is included in $(A \cup B \cup X) - B(x_1, y_2)$, in contradiction to the minimality of $X$ if $i=2$.
Therefore $i = 1$.

By Theorem~\ref{haupt} we have $G = G[A \cup B' \cup W \cup Z]$, in
contradiction to the fact that 
    $$B(x_1, y_2) \cap (A \cup B' \cup W \cup Z) = \emptyset.$$ \qed

\begin{lem}
\label{next}
There is no directed $\overline{A \cup B}$-arc joining a vertex in an $A \overline{B}$-arc to a vertex in a $B \overline{A}$-arc.
\end{lem}

\begin{figure}
\begin{center}
\setlength{\unitlength}{1cm}
\begin{picture}(14.5,7)
\put(0,0)
{\scalebox{0.30}{\includegraphics{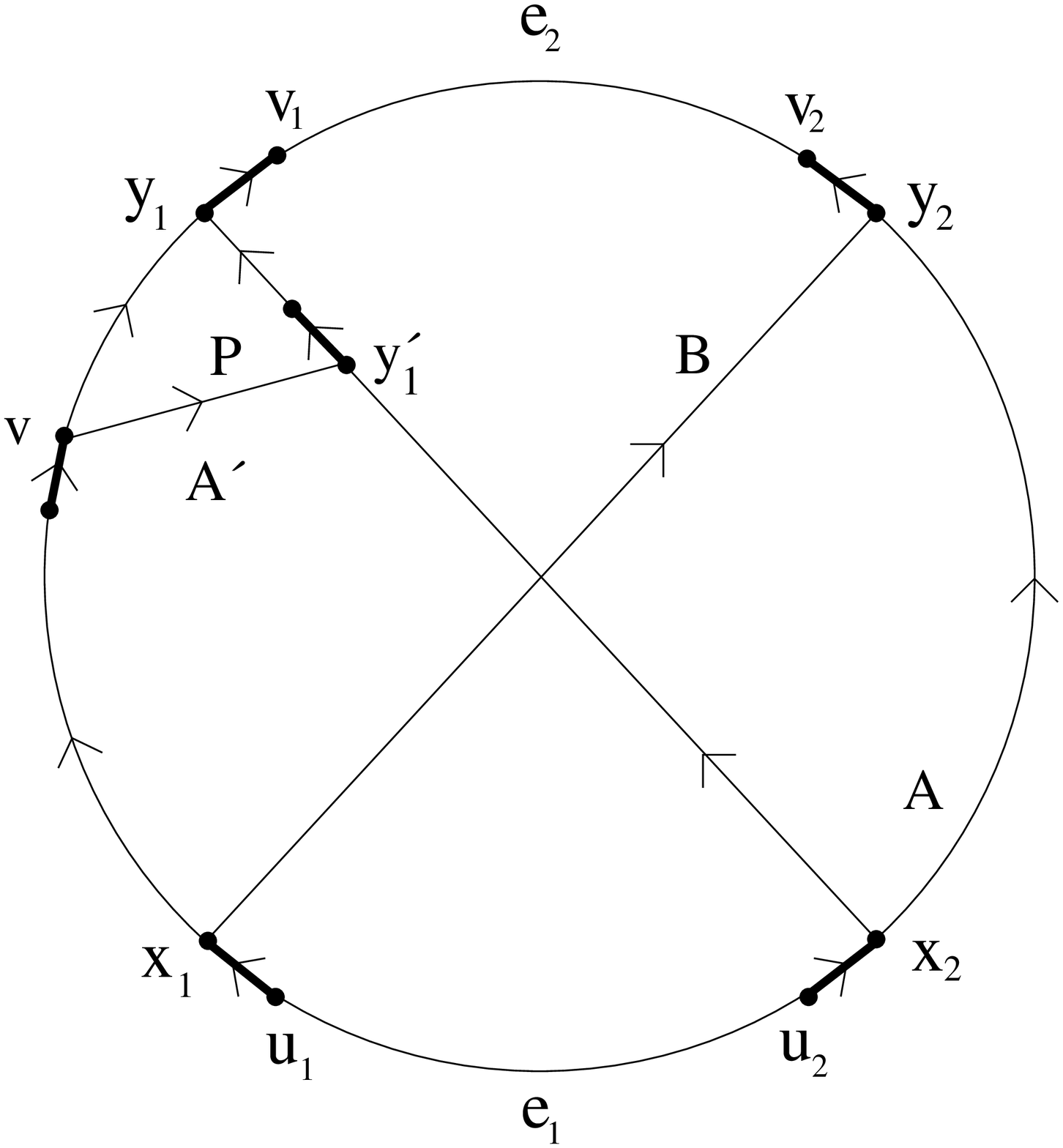}}}
\put(8,0)
{\scalebox{0.30}{\includegraphics{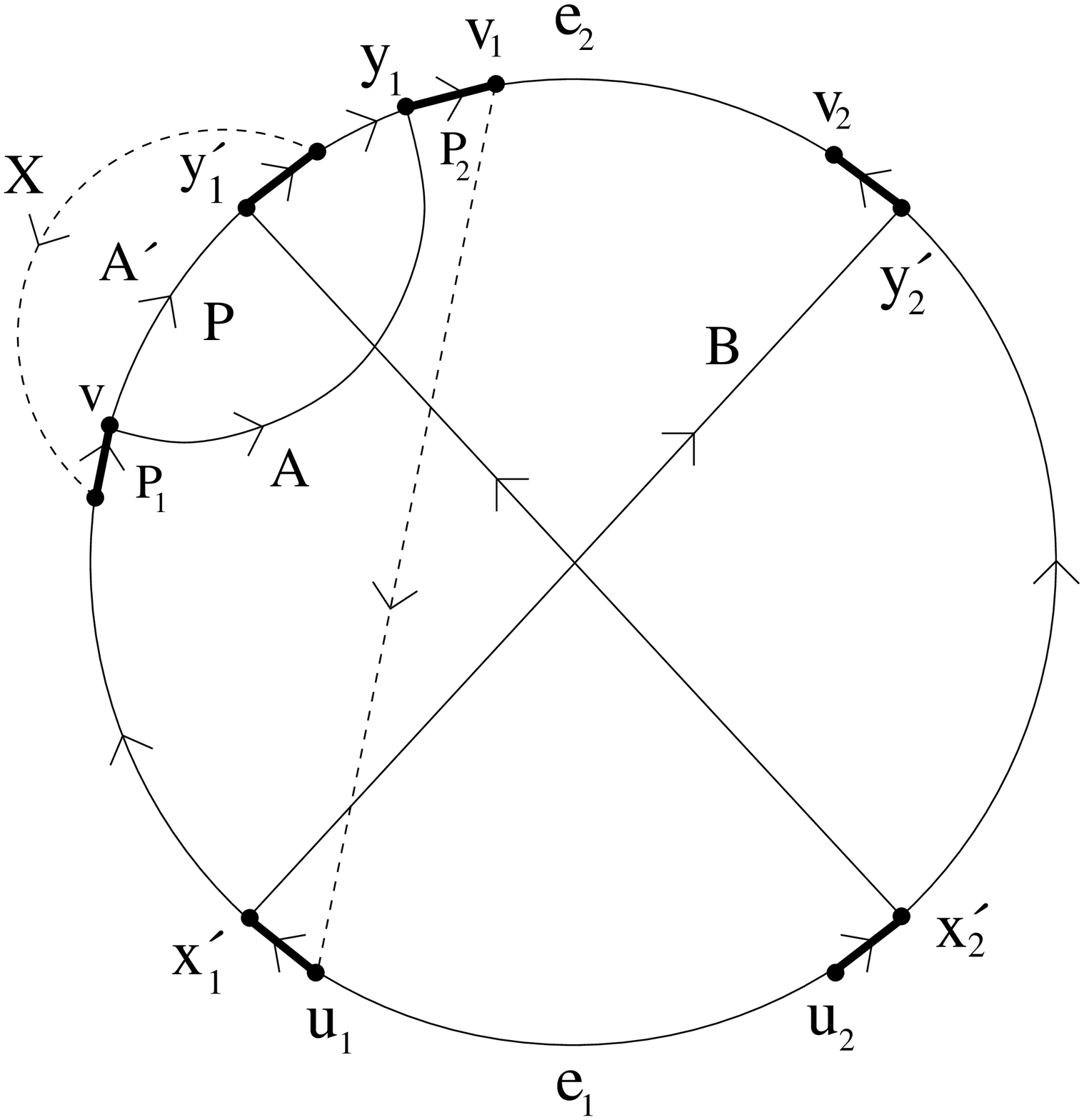}}}
\end{picture}
\end{center}
\caption{$P$ is a directed arc joining a vertex in  $A(x_1,y_1)$ to a vertex in $B(x_2,y_1)$.}
\label{AbBtoBbA}
\end{figure}

{\it Proof}: By symmetry it suffices to prove that no $\overline{A \cup B}$-arc is directed from a vertex in $A(x_1, y_1)$ 
to a vertex in $B(x_2, y_1)$. Suppose such an arc $P$ exists, joining a vertex $v \in VA(x_1, y_1)$ to a vertex $y_1' \in VB(x_2, y_1)$. 
(See Figure~\ref{AbBtoBbA}.)
Let
  $$ A' = A(u_1, v) \cup P \cup B(y_1', v_1) \cup \{e_2\} \cup A(u_2, v_2) \cup \{e_1\}. $$
This is an $f$-alternating circuit containing $e_1$ and $e_2$. Observe that $A(v, y_1)$ is the only $A \overline{A'}$-arc. It follows
that $A$ and $A'$ have the same clockwise parity.  Henceforth $A'$ will play the r\^ole previously assumed by $A$.
(See Figure~\ref{AbBtoBbA},
second picture, where $A'$ is drawn as a circle.) The circuit
$B$ will play the same r\^ole as before, but we define $x_1' = x_1$, $x_2' = x_2$ and $y_2' = y_2$. Note that $A(v, y_1)$ is an
$\overline{A' \cup B}$-arc. Therefore, by Lemma~\ref{include}, for some $i, j \in \{1,2\}$
there exists a directed minimal path $X$ from $y_i'$ to $x_j'$ including $A(v, y_1)$. 
Let $P_1$ be the $A'X$-arc with terminus $v$ and $P_2$ be the $A'X$-arc with origin 
$y_1$. Then $P_1 <_X P_2$ and $P_1 <_{A'(u_1,v_1)} P_2$ in contradiction to the 
minimality of $X$. \qed

\begin{lem}
\label{afternext}
For each $i, j \in \{1,2\}$ there is no directed $\overline{A \cup B}$-arc from a vertex in $A(y_i, v_i)$ to a vertex in $A(u_j, x_j)$.
\end{lem}

{\it Proof}: It suffices to prove the lemma for $i = j = 1$. Suppose such a directed arc $P$ exists. Let $P$ be directed from vertex $v$ to vertex $u$. Define
    $$C = A(u, v) \cup P. $$
Let $f' = f + C$,
\begin{eqnarray*}
  A' & = & A + C \\
     & = & P \cup A(v,v_1) \cup \{e_2\} \cup A(u_2, v_2) \cup \{e_1\} \cup A(u_1,u)  
\end{eqnarray*}
and 
\begin{eqnarray*}
  B' & = & B + C \\
     & = & P \cup B(v,v_1) \cup \{e_2\} \cup B(x_1, v_2) \cup A(x_1, y_1) \cup B(u_2, y_1) \cup \{e_1\} \cup B(u_1,u).
\end{eqnarray*}
Thus $A'$ and $B'$ are $f'$-alternating circuits of opposite clockwise parity containing $e_1$ 
and $e_2$. However $A(x_2, y_2)$ is the only 
$A' \overline{B'}$-arc. This result contradicts the assumption that for every choice of $A$ and $B$ there are at least two $AB$-arcs. \qed

\begin{figure}
\begin{center}
\setlength{\unitlength}{1cm}
\begin{picture}(14.5,7)
\put(0,0)
{\scalebox{0.30}{\includegraphics{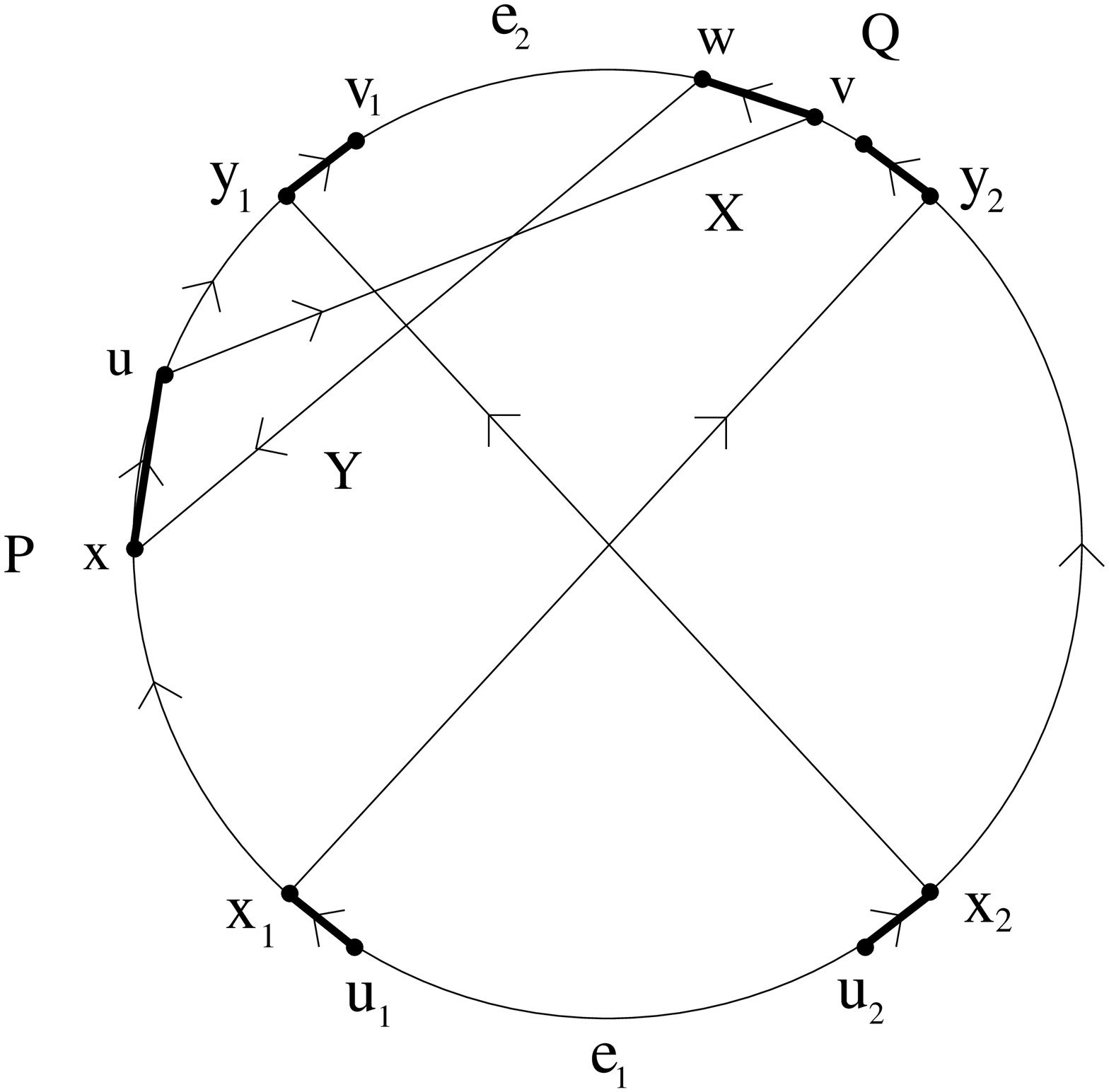}}}
\put(8,0)
{\scalebox{0.30}{\includegraphics{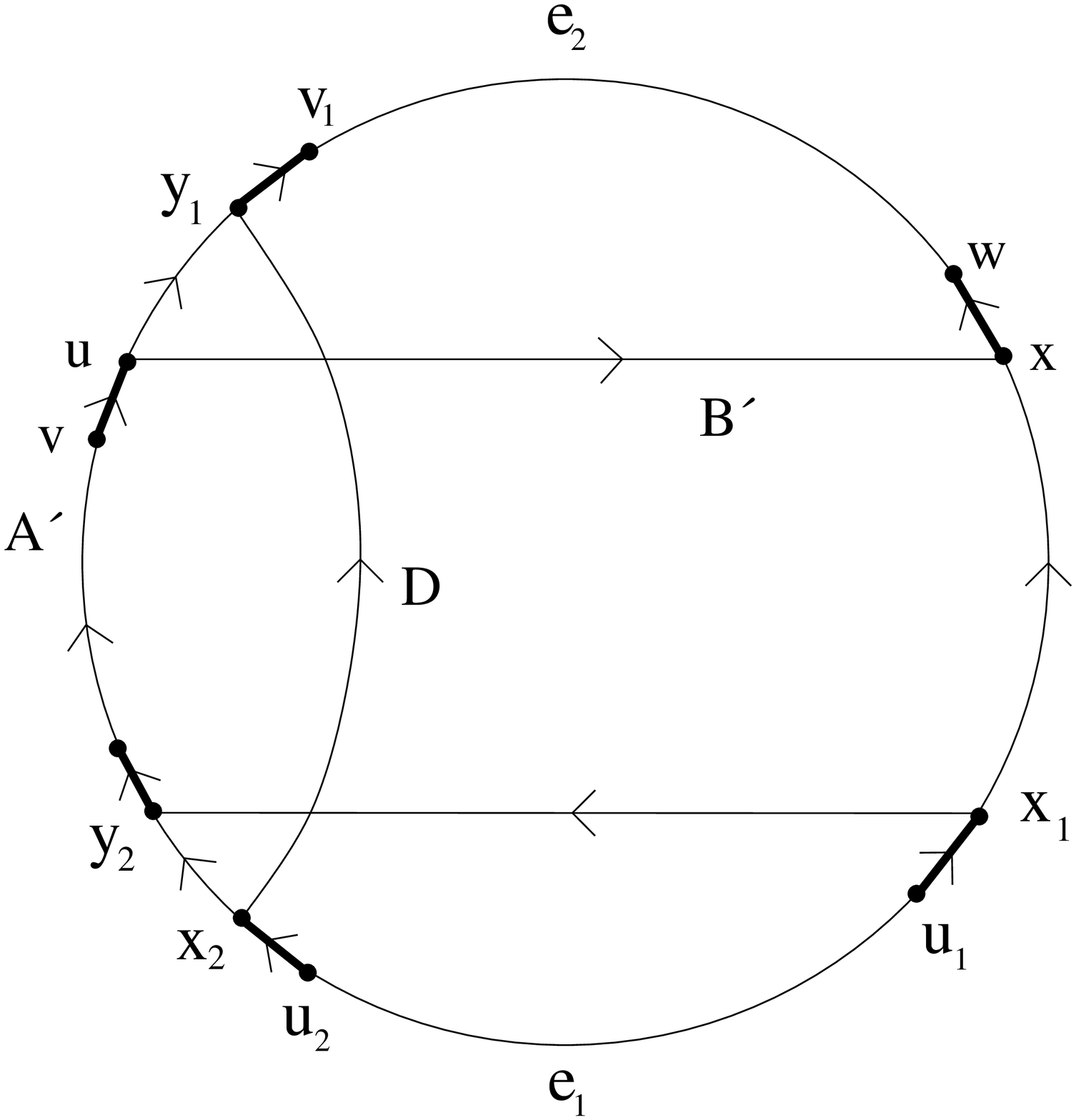}}}
\end{picture}
\end{center}
\caption{The situation in Lemma~\ref{after3} before and after reorienting $C$.}
\label{020}
\end{figure}

\begin{lem}
\label{after3}
Let $P$ be a directed $A \overline{B}$-arc or a directed $B \overline{A}$-arc and let $Q$ be a directed $AB$-arc in $H$ having neither end in common with an end of $P$.
Then there is no pair $X, Y$ of $\overline{A \cup B}$-arcs such that $X$ is directed from a vertex $u \in VP$ to a vertex $v \in VQ$, $Y$ is directed from a vertex $w \in VQ$
to a vertex $x \in VP$, $x <_P u$ and $v <_Q w$. (See Figure~\ref{020}.) 
\end{lem}

{\it Proof}: By the symmetry between $A$ and $B$ we may assume that $P = A(x_1, y_1)$.  Therefore $Q = A(u_2, x_2)$ or $Q = A(y_2, v_2)$.
By symmetry we may assume the latter case obtains.  

Suppose $X$ and $Y$ exist. Let 
   $$C = A(v,w) \cup Y \cup A(x,u) \cup X. $$
Define $f' = f + C$, 
\begin{eqnarray*}
  A' & = & A + C \\ 
     & = & A(w, v_2) \cup \{e_2\} \cup A(u, v_1) \cup X \cup A(u_2, v) \cup \{e_1\} \cup A(u_1, x) \cup Y
\end{eqnarray*}
and 
\begin{eqnarray*}
  B' & = & B + C \\
     & = & B(w, v_2) \cup \{e_2\} \cup B(u_2, v_1) \cup \{e_1\} \cup B(u_1, v) \cup X \cup A(x, u) \cup Y.
\end{eqnarray*}
Then $A'$ and $B'$ are $f'$-alternating circuits containing $e_1$ and $e_2$ and having 
opposite clockwise parity. Reorient $C$ and define
  $$ D = B'(u_2, v_1) \cup \{e_2\} \cup A'(u_1, v_2) \cup \{e_1\}. $$
(See Figure~\ref{020}.)
Then $D$ is $f'$-alternating and contains $e_1$ and $e_2$. If $D$ is of opposite clockwise parity to $A'$ then we have a contradiction because there is a unique
$A' \overline{D}$-arc $A'(x_2, y_1)$; otherwise $B'$ and $D$ have opposite clockwise parity and there is another contradiction since $B'(x_1, x)$ is the only 
$B' \overline{D}$-arc.  \qed

\section{Production of a list of cases to consider}

We now introduce a notation to describe a minimal directed path $X$ from $y_i$ to $x_j$ 
for $i, j \in \{1,2\}$. Traversed from $y_i$ to $x_j$, $X$ meets a succession of 
directed $(A \cup B)X$-arcs in $H$.
The {\it trace} of $X$ is the sequence obtained from $X$ by recording:
   $0$ for each $A(x_j, y_i)X$-arc,
   $0'$ for each $A(x_{3-j}, y_{3-i})X$-arc, 
   $1$ for each $B(x_{3-j}, y_i)X$-arc, 
   $1'$ for each $B(x_j, y_{3-i})X$-arc, 
   $2$ for each $(A \cap B)(y_{3-i}, v_{3-i})X$-arc, 
   $2'$ for each $(A \cap B)(u_{3-j}, x_{3-j})X$-arc.

\begin{figure}
\begin{center}
\setlength{\unitlength}{1cm}
\mbox{\scalebox{0.40}{%
 \includegraphics{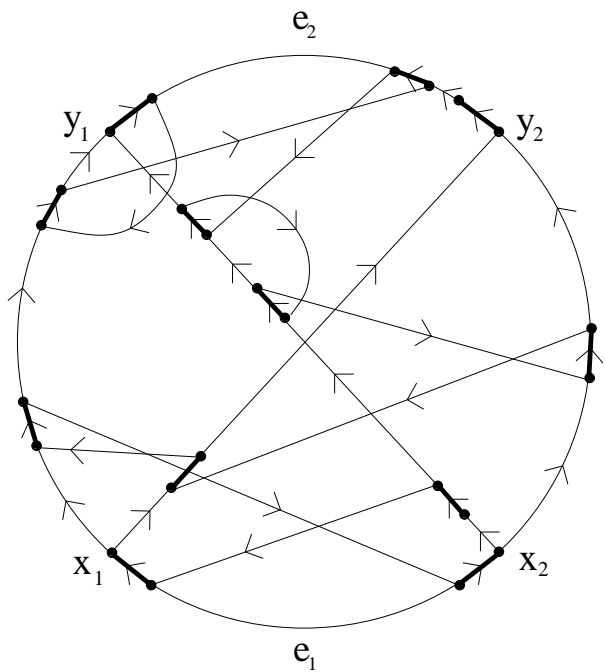}}}%
\end{center}
\caption{A directed minimal path with trace $02110'1'02'1$.}
\label{traceex}
\end{figure}

Figure~\ref{traceex} shows an example of a directed minimal path from $y_1$ to $x_1$ with trace $0 2 1 1 0' 1' 0 2' 1$.
By Lemma~\ref{minimal}(a) the graph $G[A \cup B \cup X]$ is determined up to homeomorphism by the trace of $X$. In particular, there are
a unique $A(y_i, v_i)X$-arc and a unique $A(u_j, x_j)X$-arc. 

\begin{lem}
\label{1to0}
Let $W$ be a string over $\{0,0',1,1',2,2'\}$.\\ 
(a) It is possible to choose $A$, $B$, $f$, $x_1$, $x_2$, $y_1$, $y_2$ and a directed minimal path $X$ from $y_1$ to $x_1$ in $G$ 
such that the trace of $X$ is $0 W$ if and only if it is possible, without altering $A \cup B \cup X$, to choose   
$A$, $B$, $f$, $x_1$, $x_2$, $y_1$, $y_2$ and a 
directed minimal path $X$ from $y_1$ to $x_1$ in $G$ such that the trace of $X$ is $1 W$.\\
(b) It is possible to choose $A$, $B$, $f$, $x_1$, $x_2$, $y_1$, $y_2$  and a directed minimal path $X$ from $y_1$ to $x_1$ 
in $G$ such that the trace of $X$ is $W 0$ if and only if it is possible, without altering $A \cup B \cup X$, to choose 
$A$, $B$, $f$, $x_1$, $x_2$, $y_1$, $y_2$  and a 
directed minimal path $X$ from $y_1$ to $x_1$ in $G$ such that the trace of $X$ is $W 1'$. \\
\end{lem}

\begin{figure}
\begin{center}
\setlength{\unitlength}{1cm}
\begin{picture}(14.5,7)
\put(0,0)
{\scalebox{0.30}{\includegraphics{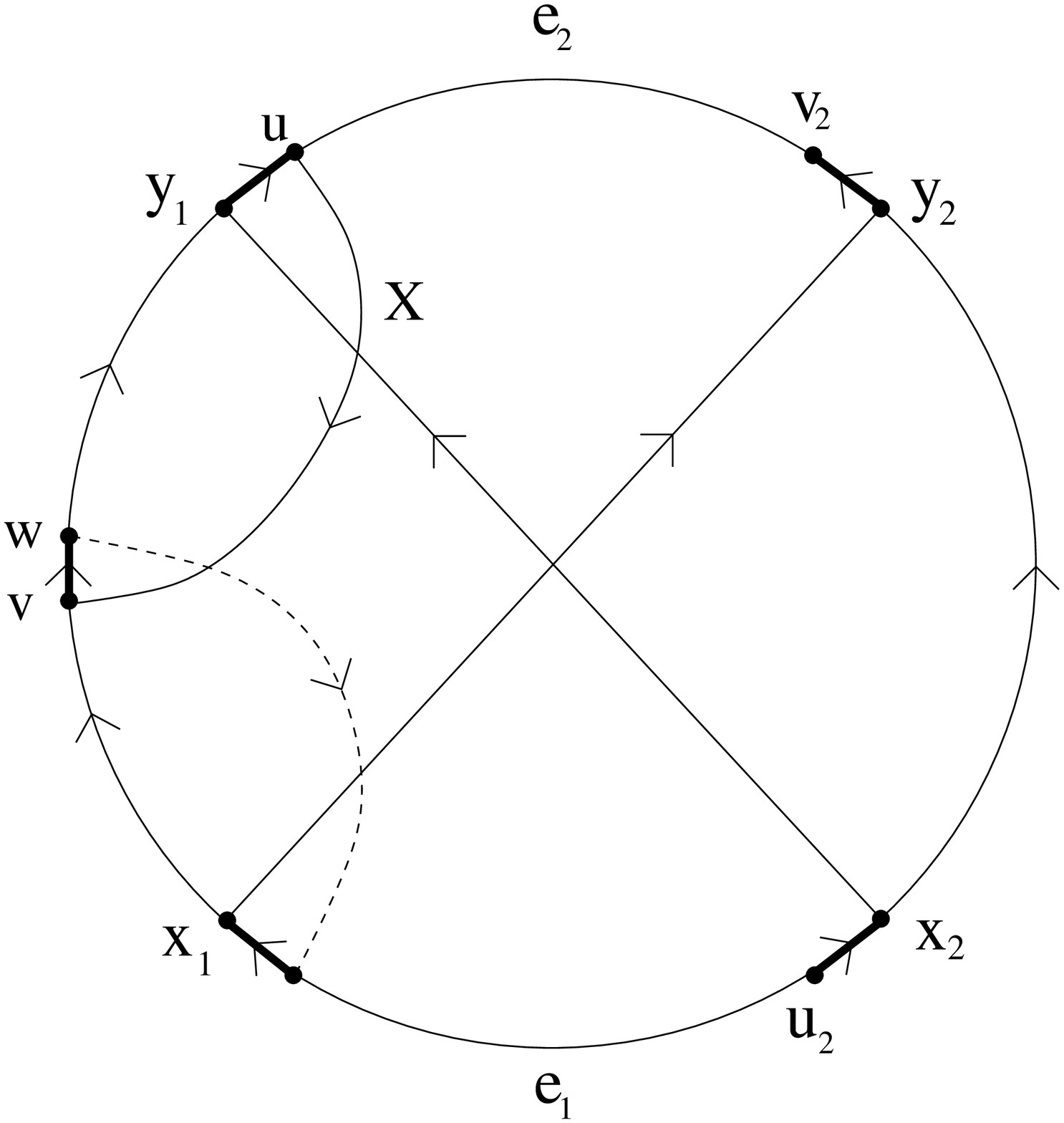}}}
\put(8,0)
{\scalebox{0.30}{\includegraphics{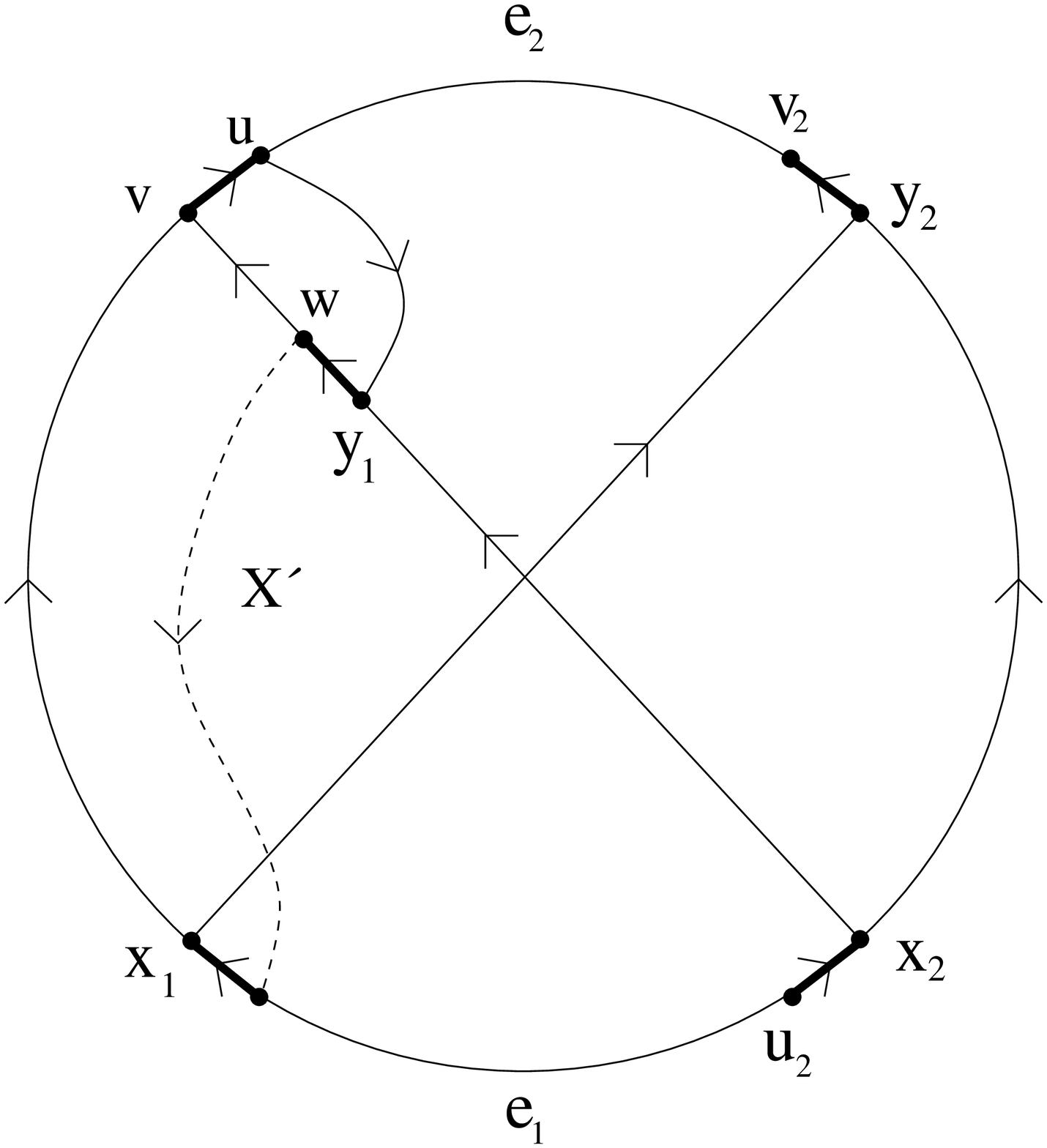}}}
\end{picture}
\end{center}
\caption{The situation in Lemma~\ref{1to0} before and after reorienting $C$.}
\label{0to1}
\end{figure}

{\it Proof}: By symmetry, it suffices to prove (a). Suppose the trace of $X$ is $0 W$.
There is an $A(x_1,y_1)X$-arc that corresponds to the first $0$ in the trace of 
$X$. Let $v$ be its origin, and let 
$$
C = X(y_1,v) \cup A(v,y_1) .
$$
Let $u$ be the terminus of the $A(y_1, v_1)X$-arc and $w$ the terminus of the 
$A(x_1, y_1)X$-arc with origin $v$. (See Figure~\ref{0to1}.) Define $f' = f + C$,
\begin{eqnarray*}
  A' & = & A + C \\
     & = & C(u, v) \cup A(u_1, v) \cup \{e_1\} \cup A(u_2, v_2) \cup \{e_2\} \cup A(u, v_1),
\end{eqnarray*}
\begin{eqnarray*}
  B' & = & B + C \\
     & = & C(u,y_1) \cup B(u_2,y_1) \cup \{e_1\} \cup B(u_1, v_2) \cup \{e_2\} 
     \cup B(u, v_1)
\end{eqnarray*}
and
\begin{eqnarray*}
  X' & = & X + A(v, y_1) \\
     & = & C(w,v) \cup X(w, x_1).
\end{eqnarray*}
Then $A'$ and $B'$ are $f'$-alternating circuits containing $e_1$ and $e_2$ and having opposite 
clockwise parity. Moreover there are exactly two $A'B'$-arcs, one containing $e_1$ and 
the other containing $e_2$, and  
the vertices of degree 3 in $G[A' \cup B']$ are $x_1$, $x_2$, $v$, $y_2$.
After reorientation of $C$, $X'$ is a directed path from $v$ to $x_1$. 
The trace of $X'$ is $1 W$, as $B'(y_1, w)$
is a $B'(x_2,v)X'$-arc which replaces the 
$A(x_1,y_1)X$-arc $A(v,w)$. (See Figure~\ref{0to1}.)
Moreover $X'$ is minimal: $X'$ satisfies the condition in 
Lemma~\ref{minimal}(b) since $X$ does and 
$$X'(v,w) \cap (A' \cup B') = C(v,u) \cup C(y_1,w).$$ 

For the converse in (a) note that $f=f'+C$, $A=A'+C$, $B=B'+C$ and $X=X'+A(v,y_1)$. \qed

\begin{figure}
\begin{center}
\setlength{\unitlength}{1cm}
\begin{picture}(14.5,7)
\put(0,0)
{\scalebox{0.30}{\includegraphics{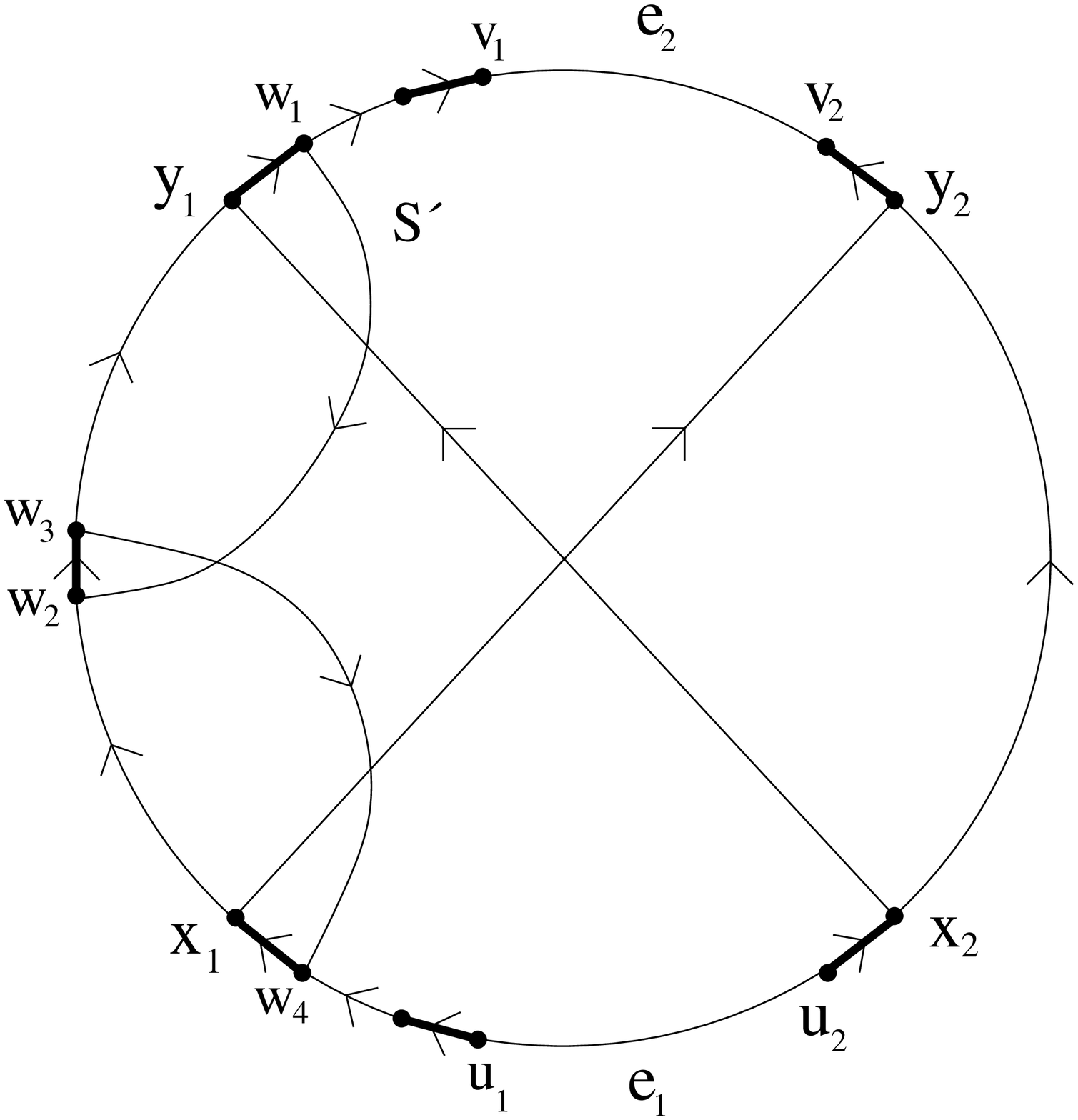}}}
\put(8,0)
{\scalebox{0.30}{\includegraphics{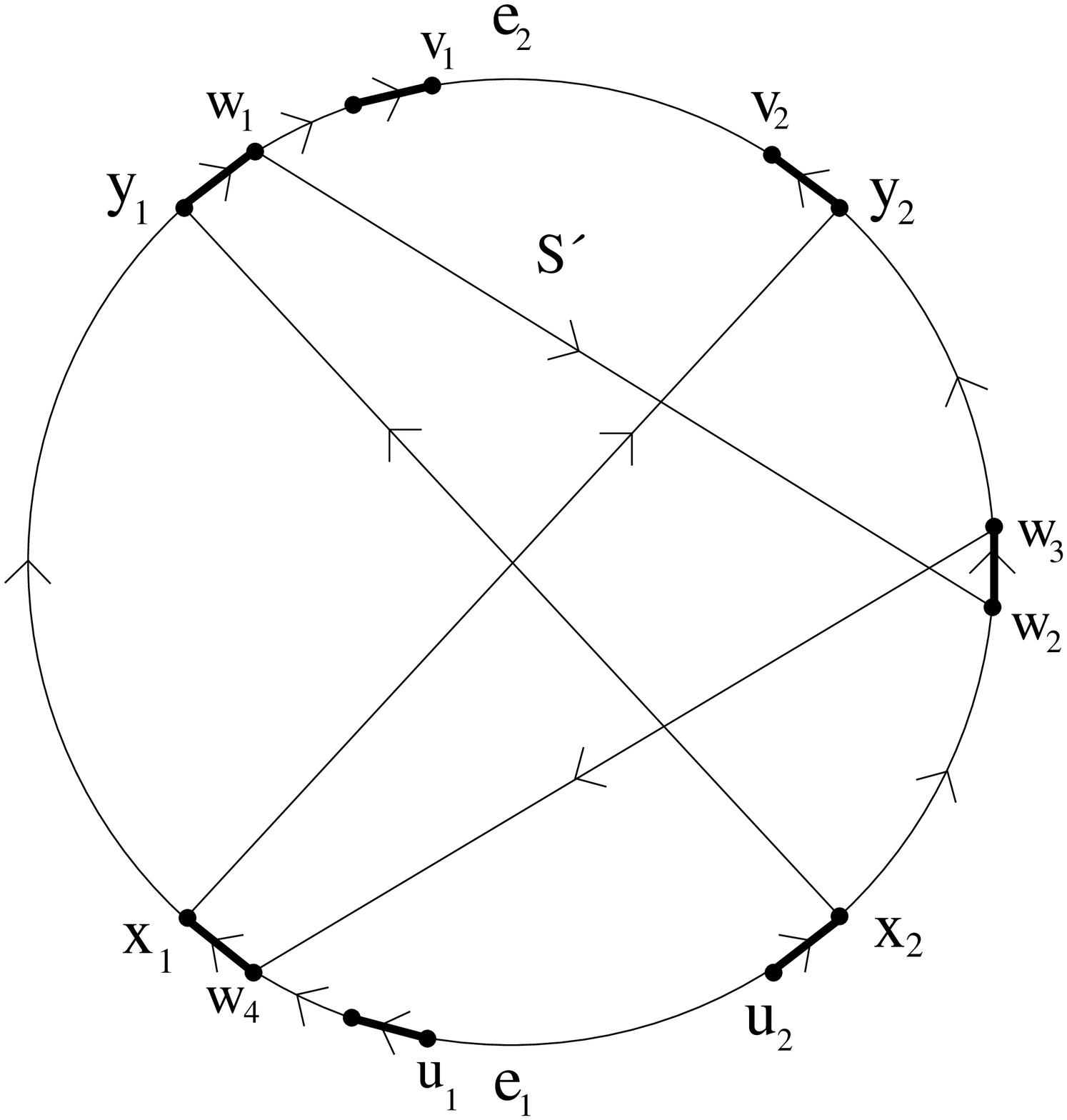}}}
\end{picture}
\end{center}
\caption{A directed minimal path $S'$ with trace $0$ and with trace $0'$.}
\label{S0}
\end{figure}

\begin{lem}
For some choice of $A$, $B$, $x_1$, $x_2$, $y_1$, $y_2$, $f$ there is a directed minimal path $S'$ from $y_1$ to $x_1$ with trace $0$ or $0'$. 
(See Figure~\ref{S0}.)
\end{lem}

{\it Proof}: We choose $A$, $B$, $x_1$, $x_2$, $y_1$, $y_2$, $f$ and a directed path $S'$ from $y_1$ to $x_1$  so that $A \cup B \cup S'$ is minimal. Thus $S'$ is minimal.

Suppose the trace of $S'$ contains $2$. There is an $(A \cap B)(y_2, v_2)S'$-arc; let $v$ be its terminus. Included in the set
$A(y_2, v) \cup S'(v, x_1)$ is a directed path from $y_2$ to $x_1$, in contradiction to the minimality of $A \cup B \cup S'$ since
$S'(y_1, v) - (A \cup B) \neq \emptyset$. Therefore the trace of $S'$ does not contain $2$, and similarly does not contain $2'$.

The trace of $S'$ contains none of $00$, $11$, $0'0'$, $1'1'$ by Lemma~\ref{long}, none of $00'$, $0'0$, $11'$, $1'1$ by Lemma~\ref{afterlong},
none of $01$, $01'$, $0'1$, $0'1'$, $10$, $10'$, $1'0$, $1'0'$ by Lemma~\ref{next}, and is non-empty by Lemma~\ref{afternext}.
We infer that the trace of $S$ is one of $0$, $0'$, $1$, $1'$. By Lemma~\ref{1to0} 
the case that the trace of $S'$ is $1$ or $1'$ can be reduced to the case 
that the trace of $S'$ is $0$.  \qed

\smallskip 

Because of this lemma we henceforth assume that the trace of $S$ is $0$ or 
$0'$. Given this choice for the trace of $S$, we now turn our attention 
to the trace of $T$.

In the following we produce a finite list of possible traces of $T$ 
and therefore a finite list of graphs we will consider in the following section.
In order to do so we distinguish the two cases
$S \cap T = \emptyset$ and $S \cap T \not= \emptyset$. First we 
assume that $S \cap T = \emptyset$. 

\begin{lem}
\label{1to0v}
Let $W$ be a string over $\{0,0',1,1',2,2'\}$, and let $S \cap T = \emptyset$.\\ 
(a) Suppose that the trace of $T$ is $1 W$. Then there exist $A$, $B$, $f$, $x_1$, $x_2$, $y_1$, $y_2$, a directed minimal path 
from $S'$ $y_1$ to $x_1$ with trace $0$ or $0'$ and a directed minimal path $T'$ from
$y_2$ to $x_2$ with trace $0W$. \\
(b) Suppose  that the trace of $T$ is $W 1'$. Then there exist $A$, $B$, $f$, $x_1$, $x_2$, $y_1$, $y_2$, a directed minimal path 
$S'$ from $y_1$ to $x_1$ with trace $0$ or $0'$ and a directed minimal path  $T'$ from
$y_2$ to $x_2$ with trace $W 0$.
\end{lem}

{\it Proof}: That $T'$ exists follows by a proof similar to that of 
the corresponding assertion in Lemma~\ref{1to0}. The reorientation of the corresponding $f$-alternating circuit 
$C$ does not affect $S$, since $S \cap C = \emptyset$. Therefore we may take $S' = S$. \qed

\smallskip

Thus we can assume that the first symbol in the trace of $T$ is not $1$ and 
that the last symbol in the trace of $T$ is not $1'$. In the following we will refer to 
this property of $T$ as $(A)$.

\smallskip

\begin{lem}
\label{before}
Suppose $S \cap T = \emptyset$ and the trace of $T$ contains one of $20'$, $21'$, $12'$, $0'2'$. Then there 
exist $A'$, $B'$, $f'$, a directed minimal path $S'$ from $y_1$ to $x_1$ with 
trace $0$ or $0'$  and a directed minimal path $T'$ from $y_2$ to $x_2$ such that 
$S' \cap T' \not= \emptyset$. 
\end{lem}

\begin{figure}
\begin{center}
\setlength{\unitlength}{1cm}
\begin{picture}(14.5,7)
\put(0,0)
{\scalebox{0.30}{\includegraphics{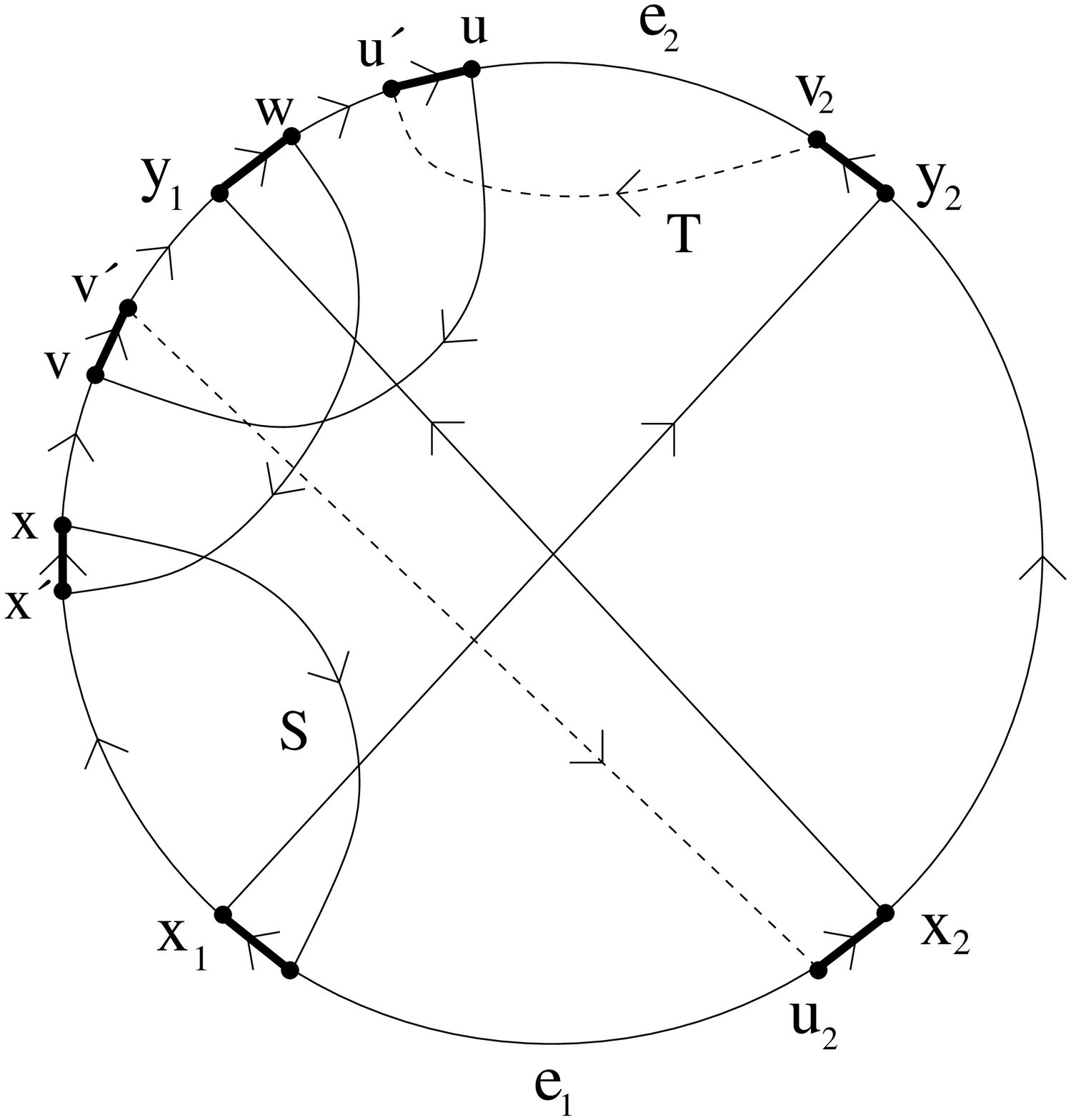}}}
\put(8,0)
{\scalebox{0.30}{\includegraphics{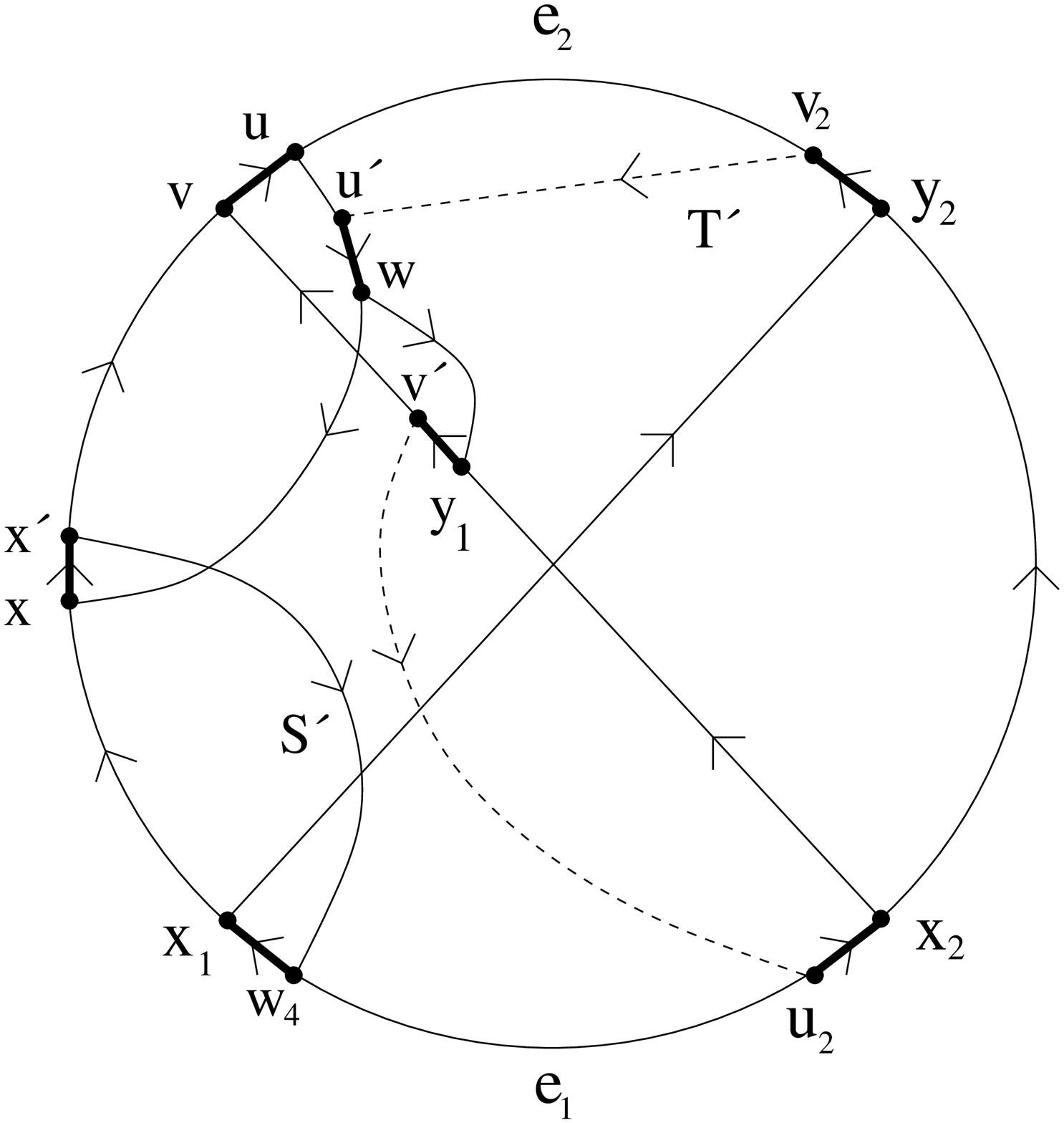}}}
\end{picture}
\end{center}
\caption{The situation in Lemma~\ref{before} before and after reorienting $C$.}
\label{T20z}
\end{figure}

{\it Proof}: It suffices to consider the case where the trace of $T$ contains $20'$, for the other cases are similar. 
In this case there is a $T \overline{A \cup B}$-arc with origin 
$u \in VA(y_1, v_1)$ and terminus $v \in VA(x_1, y_1)$. Let $u'$ be the origin of the $A(y_1, v_1)T$-arc with terminus $u$, and $v'$ the terminus 
of the $A(x_1, y_1)T$-arc with origin $v$. 
Let $w$ be the terminus of the unique $A(y_1, v_1)S$-arc. (See Figure~\ref{T20z}.) 

Since $S \cap T = \emptyset$ we have $w <_{A(y_1, v_1)} u'$. If $S$ has trace $0$, then let
$x$ be the terminus of the unique $A(x_1, y_1)S$-arc and $x'$ its origin. We have $x <_{A(x_1, y_1)} v$: otherwise if we define 
  $$ S^* = A(y_1, u) \cup T(u, v) \cup A(v, x) \cup S(x, x_1) $$
then $G = G[A \cup B \cup S^* \cup T]$ by Theorem~\ref{haupt}, in contradiction to the fact that 
  $$ S(w, x') \cap (A \cup B \cup S^* \cup T) = \emptyset. $$
In any case, define
  $$ C = T(u, v) \cup A(v, u). $$
This is an $f$-alternating circuit such that $C \cap S = A(y_1, w)$. Reorient $C$ and  define $f' = f + C$,
\begin{eqnarray*}
  A' & = & A + C \\
     & = & C(v,u) \cup A(u_1, v) \cup \{e_1\} \cup A(u_2, v_2) \cup \{e_2\} \cup A(u, v_1),
\end{eqnarray*}
\begin{eqnarray*}
  B' & = & B + C \\
     & = & C(y_1,u) \cup B(u_2, y_1) \cup \{e_1\} \cup B(u_1, v_2) \cup \{e_2\} \cup B(u, v_1),
\end{eqnarray*}
\begin{eqnarray*}
  S' & = & S + C(v,y_1) \\
     & = & C(v,w) \cup S(w, x_1)
\end{eqnarray*}
and 
\begin{eqnarray*}
  T' & = & T + C \\
     & = & T(y_2, u') \cup C(u',v') \cup T(v', x_2).
\end{eqnarray*}
Then $A'$ and $B'$ are $f'$-alternating circuits containing $e_1$ and $e_2$ and having opposite clockwise parity.
There are exactly two $A'B'$-arcs, the vertices of degree $3$ in $G[A' \cup B']$ 
are $x_1$, $x_2$, $v$ and $y_2$, 
$S'$ is a directed path from $v$ to $x_1$ and $T$ is a directed path from $y_2$ to $x_2$ (see Figure~\ref{T20z}).
Moreover $S$ and $S'$ have equal trace  and  $S' \cap T' = C(u', w) \neq \emptyset$. 
Finally $S'$ and $T'$ are minimal: both satisfy the condition in Lemma~\ref{minimal}(b)
since $S'(v,w) \cap (A' \cup B')=C(v,u)$ and $T'(u',v') \cap (A' \cup B') = C(y_1,v')$. \qed

\smallskip

Thus we assume that the trace of $T$ contains none of $20'$, $21'$, $12'$, $0'2'$, if $S \cap T = \emptyset$.
In the following we will refer to this property of $T$ as $(B)$. 

The following lemma gives a complete list of graphs to be considered if 
$S \cap T = \emptyset$. We use $*$ to denote an arbitrary string of symbols, and $\Lambda$ to denote the empty string.

\begin{lem}
\label{nointersection}
Suppose $S \cap T = \emptyset$ and that $T$ has properties $(A)$ and $(B)$. Then the trace of $T$ is one of $0$, $0'$, $20$, $21$, $02'$, $1'2'$.
\end{lem}   

{\it Proof}: First we see that the symbols in the trace of $T$ alternate between the sets $\{0, 0', 1, 1'\}$ and $\{2, 2'\}$, for the
trace of $T$ contains none of $00$, $11$, $22$, $0'0'$, $1'1'$, $2'2'$ by Lemma~\ref{long}, 
none of $00'$, $0'0$, $11'$, $1'1$ by Lemma~\ref{afterlong},
none of $01$, $01'$, $0'1$, $0'1'$, $10$, $10'$, $1'0$, $1'0'$ by Lemma~\ref{next}, 
and does not contain $22'$ by Lemma~\ref{afternext} or $2'2$ by Lemma~\ref{minimal}(a).

Next we show that the trace of $T$ does not contain both $2$ and $2'$.
Suppose that the trace of $T$ contains $2'*2$. Thus $x_1,y_1 \in VT$ 
by Lemma~\ref{minimal}(a). This is a
contradiction to $S \cap T = \emptyset$.
Suppose that the trace of $T$ contains $2*2'$. Choose $2$ and $2'$ in the trace of $T$, 
with the chosen $2'$ appearing later than the chosen $2$.
Let $u$ be the terminus of the 
$A(y_1,v_1)T$-arc that corresponds to the chosen $2$ and let 
$v$ be the origin of the $A(u_1,x_1)T$-arc that corresponds to the chosen 
$2'$. There is a directed path from $y_1$ to 
$x_1$ included in 
$$ A(y_1,u) \cup T(u,v) \cup A(v,x_1). $$ 
By Theorem~\ref{haupt} we have $G=G[A \cup B \cup T]$, in contradiction to 
$$ (S - (A \cup B)) \cap (A \cup B \cup T) = \emptyset.$$

We distinguish the following cases:
\begin{enumerate}
\item The trace of $T$ contains neither $2$ nor $2'$.
\item The trace of $T$ contains $2$ and consequently does not contain $2'$.
\item The trace of $T$ contains $2'$ and consequently does not contain $2$.
\end{enumerate}

{\it Case 1}: By Lemma~\ref{afternext} the trace is not empty and therefore the
trace of $T$ is one of $0$, $0'$, $1$, $1'$ in this case. By $(A)$ the case that the trace is 
$1$ or $1'$ is not possible.

{\it Case 2}: The symbols in the trace of $T$ alternate between the sets $\{2\}$ and $\{0, 0', 1, 1'\}$. 
By $(B)$ and Lemma~\ref{afternext} every $2$ in the trace must be immediately followed 
by $0$ or $1$. Therefore the trace contains at most one $2$ by Lemma~\ref{after3} 
and consequently exactly one $2$. In fact the trace of $T$ is $x 2 y$, 
where $x \in \{\Lambda,0,0',1,1'\}$ and $y \in \{0,1\}$. 

We show that $x=\Lambda$. We have $x \notin \{0',1'\}$, for otherwise 
$y_1 \in VS \cap VT$ by Lemma~\ref{minimal}(a). Suppose $x =0$. If $y=0$, we have a contradiction 
by Lemma~\ref{after3}. If $y=1$, we have a contradiction by Lemma~\ref{after3}
also, since in this case there exist $A'$, $B'$, $f'$, $x_1'$, $x_2'$, $y_1'$, $y_2'$ in $G$ and a directed minimal 
path $T'$ from $y'_2$ to $x'_2$ with trace $121$ by Lemma~\ref{1to0}.
Similarly we obtain a contradiction if we suppose that $x=1$. 
Therefore the trace of $T$ is either $20$ or $21$, if the trace contains $2$.

{\it Case 3}: Similarly the trace of $T$ is either $02'$ or $1'2'$, if it contains $2'$. \qed

\begin{rem}
\label{symmetric1}
The case where the trace of $T$ is either $02'$ or $1'2'$ can be reduced to the case where the trace of $T$ is either 
$20$ or $21$. 
In order to see this suppose that the trace of $T$ is either $02'$ or $1'2'$ and switch to the reference orientation with respect to 
$(N,M,f)$. Then $S$ is a directed path from $y'_1=x_1$ to $x'_1=y_1$  with trace $0$ or $0'$ and $T$ is a directed path from 
$y'_2=x_2$ to $x'_2=y_2$ with trace $20$ or $21$.
\end{rem}

\smallskip

Now we consider the case where $S \cap T \not= \emptyset$.
First we prove the following consequence of Theorem~\ref{haupt}.

\begin{cor}
\label{hauptcor}
Let $S'$ be a directed path from $y_1$ to $x_1$, $T_1$ a directed path from 
$y_2$ to a vertex in $S'$ and $T_2$ a directed path from a vertex in $S'$ to 
$x_2$. Then 
$$ G=G[A \cup B \cup S' \cup T_1 \cup T_2]. $$
\end{cor}

{\it Proof}: By Theorem~\ref{haupt} we have to show that there exists a directed 
path from $y_2$ to $x_2$ in $G[A \cup B \cup S' \cup T_1 \cup T_2]$. Let 
$a$ be the terminus of $T_1$ and $b$  the origin of $T_2$. Then such a path 
is included in 
$$ T_1 \cup S(a,x_1) \cup A(x_1,y_1) \cup S(y_1,b) \cup T_2.$$ \qed

\begin{figure}
\begin{center}
\setlength{\unitlength}{1cm}
\begin{picture}(14.5,7)
\put(0,0)
{\scalebox{0.30}{\includegraphics{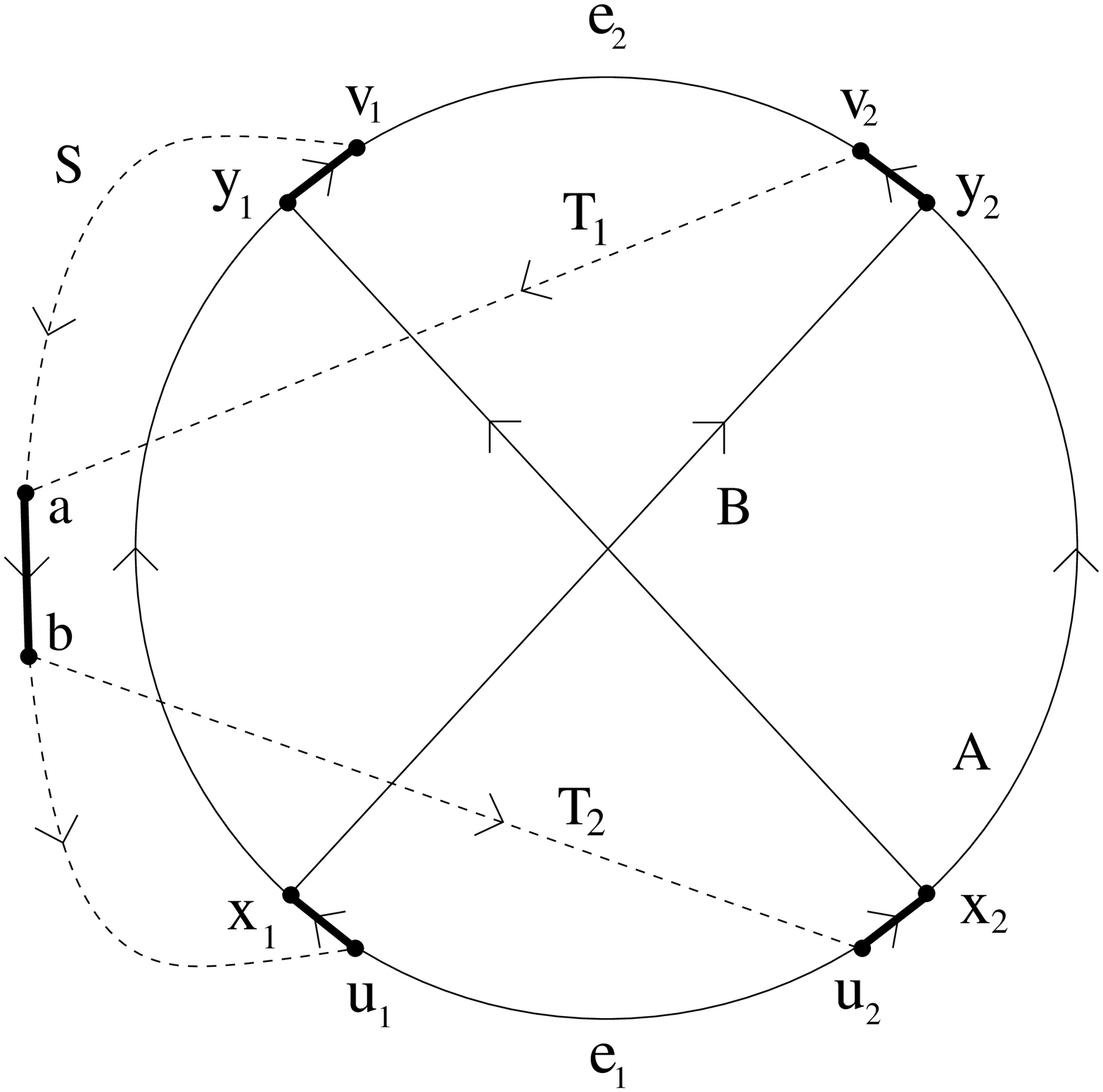}}}
\put(8,0)
{\scalebox{0.30}{\includegraphics{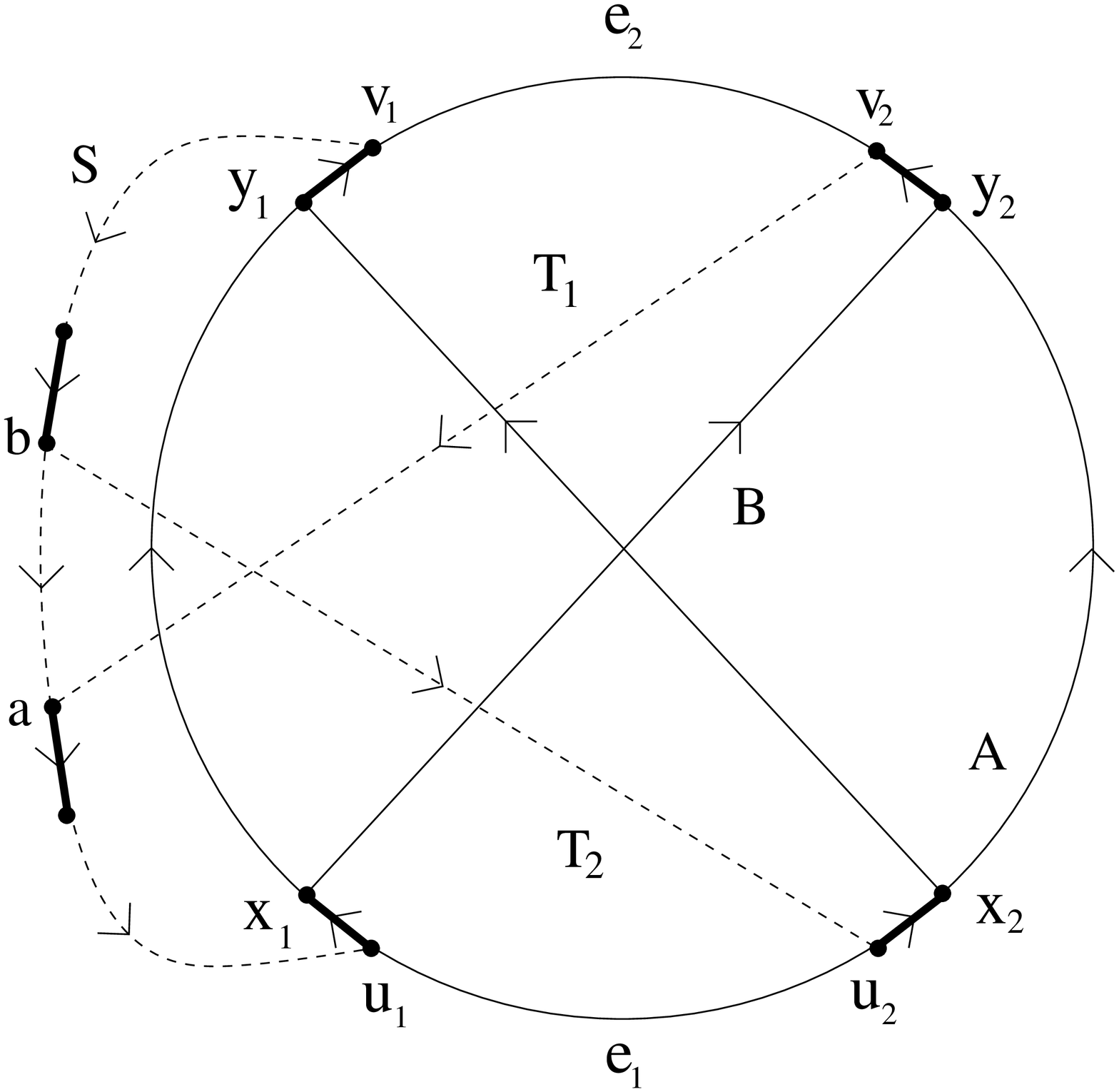}}}
\end{picture}
\end{center}
\caption{There are two possible cases: either $a <_S b$ or $b <_S a$.}
\label{hauptcorfig}
\end{figure}

Since $S \cap T \not= \emptyset$, there exists a first vertex $a$ in $T$ that 
is also in $S$, and a last vertex $b$ in $T$ that is also in $S$. (See Figure~\ref{hauptcorfig}.) 
Let $T_1=T(y_2,a)$ and $T_2=T(b,x_2)$. By Corollary~\ref{hauptcor}, 
$$ G=G[A \cup B \cup S \cup T_1 \cup T_2]. $$

We define the trace of $T_1$ and $T_2$ in a manner analogous to the definition of the trace of 
a directed minimal path from $y_2$ to $x_2$. Note that $T_1$ and 
$T_2$ satisfy the condition in Lemma~\ref{minimal}(b), since they are 
directed subpaths of the directed minimal path $T$.

\begin{lem}
\label{1to0w}
Let $W$ be a string over $\{0,0',1,1',2,2'\}$, and let $S \cap T \not= \emptyset$.\\ 
(a) Suppose that the trace of $T_1$ is $1 W$. Then there exist $A$, $B$, $f$, $x_1$, $x_2$, $y_1$, $y_2$, a directed minimal path 
from $y_1$ to $x_1$ with trace $0$ or $0'$ and a directed minimal path 
$T$ from $y_2$ to $x_2$ such that the trace of $T_1$ is $0W$. \\
(b) Suppose  that the trace of $T_2$ is $W 1'$. Then there exist $A$, $B$, $f$, $x_1$, $x_2$, $y_1$, $y_2$, a directed minimal path 
from $y_1$ to $x_1$ with trace $0$ or $0'$ and a directed minimal path 
$T$ from $y_2$ to $x_2$ such that the trace of $T_2$ is $W 0$.
\end{lem}

{\it Proof}: Similar to the proof of Lemma~\ref{1to0v}. \qed

\smallskip

Therefore we can assume that the first symbol in the trace of $T_1$ is not $1$ and that the last symbol in 
the trace of $T_2$ is not $1'$. In the following we will refer to this property of $T_1$ and $T_2$ as $(A')$.

In the directed path $S$ there are exactly $6$ vertices of degree $3$ in $G[A \cup B \cup S]$,
the first being  $y_1$ and the last being $x_1$. We label the other such vertices 
$w_1$, $w_2$, $w_3$ and $w_4$ in the order they occur when we traverse 
$S$ from $y_1$ to $x_1$. (See Figure~\ref{S0}.) 
Note that $a \notin \{w_1, w_3, x_1\}$ and $b \notin \{y_1, w_2, w_4\}$,
since the vertices of $G$ have indegree 1 or outdegree 1. 

\begin{lem}
\label{ab}
The vertices $a$ and $b$ are not both in $VS(y_1,w_1)$ and not 
both in $VS(w_4,x_1)$.
\end{lem}

{\it Proof}: By symmetry it suffices to show that 
$\{a,b\} \not\subseteq VS(y_1,w_1)$. Suppose the contrary.

First we assume that the trace of $S$ is $0$.  We define 
$$ C= S(w_1,w_2) \cup A(w_2,w_1). $$
This is an $f$-alternating circuit.
Furthermore we define $f'=f+C$, 
\begin{eqnarray*}
A' & = & A + C \\
   & = & C(w_1,w_2) \cup A(u_1,w_2) \cup \{e_1\} \cup A(u_2,v_2) \cup \{e_2\} \cup A(w_1,v_1),
\end{eqnarray*}
\begin{eqnarray*}
B' & = & B + C \\
   & = & C(w_1,y_1) \cup B(u_2,y_1) \cup \{e_1\} \cup 
         B(u_1,v_2) \cup \{e_2\} \cup B(w_1,v_1).
\end{eqnarray*}
By Lemma~\ref{parity}, $A'$ and $B'$ are $f'$-alternating circuits containing $e_1$ and $e_2$, of 
opposite clockwise parity, such that there are exactly two $A'B'$-arcs. The 
vertices of degree $3$ in $G[A' \cup B']$ are $x_1,x_2,w_2,y_2$. 

First we assume $a<_S b$. 
In this case $a=y_1$ and $b=w_1$, for otherwise we would have 
vertices of degree $2$. We 
define the paths
$$ X'=T_1 \cup C(w_3,a) \cup S(w_3,x_1) $$
and 
$$ Y'=C(b,w_2) \cup T_2.$$ 
These paths would become directed from $y_2$ to $x_1$ and from 
$w_2$ to $x_2$, respectively, if $C$ were reoriented. Therefore 
$$G=G[A' \cup B' \cup X' \cup Y']$$
by Theorem~\ref{haupt},
in contradiction to 
$$C(a,b) \cap  (A' \cup B' \cup X' \cup Y') = \emptyset.$$

Now we assume that $b <_S a$.
We reorient $C$ and define the directed path
\begin{eqnarray*}
T' & = &  T + C \\
   & = &  T_1 \cup C(a,b) \cup T_2
\end{eqnarray*}
and a directed path $S'$ from $w_2$ to $x_1$ included in 
$$ C(w_2,b) \cup T_2 \cup B(x_2,y_1) \cup C(y_1,w_3) \cup S(w_3,x_1). $$
By Theorem~\ref{haupt}
$$G=G[A' \cup B' \cup S' \cup T'],$$
in contradiction to 
$$ C(b,y_1) \cap (A' \cup B' \cup S' \cup T') = \emptyset.$$

Next we consider the case that the trace of $S$ is $0'$.
We define 
$$ D= S(w_1,w_4) \cup A(w_4,w_1). $$
This is an $f$-alternating circuit.
Furthermore we define $f''=f+D$, 
\begin{eqnarray*}
A'' & = & A + D \\
   & = & D(w_1,w_2) \cup A(u_2,w_2) \cup \{e_1\} \cup 
         A(u_1,w_4) \cup   \\ && D(w_3,w_4) \cup A(w_3,v_2) \cup \{e_2\} \cup A(w_1,v_1),
\end{eqnarray*}
\begin{eqnarray*}
B'' & = & B + D \\
   & = & D(w_1,w_4) \cup B(u_1,w_4) \cup \{e_1\} \cup 
         B(u_2,y_1) \cup \\ && D(x_1,y_1) \cup B(x_1,v_2) \cup \{e_2\} \cup B(w_1,v_1).
\end{eqnarray*}
By Lemma~\ref{parity}, $A''$ and $B''$ are $f''$-alternating circuits including $e_1$ and $e_2$, of 
opposite clockwise parity, such that there are exactly two $A''B''$-arcs. The 
vertices of degree $3$ in $G[A'' \cup B'']$ are $w_3,x_2,w_2,y_2$. 

First we assume $a <_S b$. Again we have $a=y_1$ and $b=w_1$, for otherwise $G$ would have 
vertices of degree $2$.  We 
define the paths
$$ X''= D(b,w_2) \cup T_2$$
and 
$$ Y''=T_1 \cup D(w_3,a).$$
These paths would become directed from $w_2$ to $x_2$ and from 
$y_2$ to $w_3$, respectively, if $D$ were reoriented. Therefore 
$$G=G[A'' \cup B'' \cup X'' \cup Y'']$$
by Theorem~\ref{haupt},
in contradiction to 
$$D(a,b) \cap  (A'' \cup B'' \cup X'' \cup Y'') = \emptyset. $$

Now we assume that $b <_S a$.
We reorient $D$ and define the directed path
\begin{eqnarray*}
T'' & = &  T + D \\
   & = &  T_1 \cup D(a,b) \cup T_2
\end{eqnarray*}
and a directed path $S''$ from $w_2$ to $w_3$ included in 
$$ D(w_2,b) \cup T_2 \cup B(x_2,y_1) \cup D(y_1,w_3). $$
By Theorem~\ref{haupt}
$$G=G[A'' \cup B'' \cup S'' \cup T''],$$
in contradiction to 
$$ D(b,y_1) \cap (A'' \cup B'' \cup S'' \cup T'') = \emptyset.$$  \qed

\begin{lem}
\label{trace}
The traces of $T_1$ and $T_2$ are either empty or $0$.
\end{lem}

{\it Proof}: 
By symmetry it suffices to show that the trace of $T_1$ is either empty or $0$.
First we show that the edge in $T_1$ incident on $a$ is not in $A \cup B$. This edge exists, for $T_1 \not= \emptyset$ since the trace of $S$ is 0 or $0'$.
Let $a^*$ be its origin and  suppose that $(a^*,a) \in A \cup B$. Note that $(a^*,a) \notin f$ because $S$ is $f$-alternating. Thus the 
edge of $f$ incident on $a^*$ is in $T_1$ and in $A \cup B$. Now we have 
the contradiction that $a^*$ is a vertex of degree $2$, since $a^* \notin VS \cup VT_2$. 

We use this observation to show that 
the trace of $T_1$ contains none of $0'$, $1'$, $2'$. 

Suppose that the trace of $T_1$ contains $0'$. Let $a'$ be the origin of 
an $A(x_1,y_1)T_1$-arc. Clearly $a' \not= a$.
Included in 
$$ T_1(y_2,a') \cup A(a',y_1)$$
is a directed path $T_1'$ from $y_2$ to a vertex in $S$ 
with $$T_1' \cap (T_1(a',a) - (A \cup B)) = \emptyset.$$
By Corollary~\ref{hauptcor} we have 
$$ G=G[A \cup B \cup S \cup T_1' \cup T_2]. $$
Since 
$$ (T_1(a',a) - (A \cup B)) \cap (A \cup B \cup S \cup T_1' \cup T_2) = \emptyset,$$
it follows that 
$$ T_1(a',a) - (A \cup B) = \emptyset. $$
Therefore we have the contradiction that the edge of $T_1$ incident on $a$ is in $A \cup B$. 
The proof that the trace of $T_1$ contains neither $1'$ nor $2'$ is similar. Likewise the trace of $T_2$ contains 
none of $0'$, $1$ and $2$.

Next we show that  the trace of $T_1$ does not contain 
$2$. Assume the contrary and let $v$ be the terminus of the last 
$A(y_1,v_1)T_1$-arc.
Then there is a directed path $S'$
from $y_1$ to $x_1$ included in 
$$ A(y_1,v) \cup T_1(v,a) \cup S(a,x_1) $$
with 
$$S' \cap (S(y_1,a) - (A \cup B)) = \emptyset .$$
Therefore 
$$ S(y_1,a) \subseteq A \cup B \cup T, $$
since $G=G[A \cup B \cup S' \cup T]$ by Theorem~\ref{haupt}.
Consequently, 
  $$a \notin (VS(w_1,w_2) - \{w_1\}) \cup (VS(w_3,w_4) - \{w_3\}), $$
for otherwise 
the edge of $S(y_1,a)$ incident on $a$ is in $T$ (since it is not in $A \cup B$), and we have a contradiction 
to the choice of $a$. Thus 
$$ a \in VS(y_1,w_1)  \cup (VS(w_2,w_3) - \{w_2\}) \cup (VS(w_4,x_1) - \{w_4\}). $$
Now we distinguish three cases according to which set of this union contains $a$. 

Suppose that $a \in VS(y_1,w_1)$. We already know that the trace of $T_1$ is a string over 
$\{0,1,2\}$. The symbols alternate between $2$ and members of the set $\{0,1\}$ 
for the trace of $T_1$ 
contains none of $00$, $11$, $22$ by Lemma~\ref{long} and 
neither $01$ nor $10$ by Lemma~\ref{next}. Therefore $2$ is either the last symbol or the penultimate symbol in 
the trace of $T_1$. If $2$ is the last symbol in 
the trace of $T_1$ we have a contradiction by Lemma~\ref{long}, and if 
$2$ is the penultimate symbol in the trace of $T_1$ then the last 
symbol of the trace of $T_1$ is either $0$ or $1$ and 
we have a contradiction by Lemma~\ref{minimal}(a) and Lemma~\ref{after3}.

Suppose that $a \in VS(w_4,x_1) - \{w_4\}$. Note that $a$ cannot be adjacent to $w_4$ since both vertices have indegree more than 1. Thus $b \in VS(w_4, a) - \{w_4\}$,
for otherwise there would be a vertex of degree 
$2$ in $G$. Therefore $\{a,b\} \subseteq VS(w_4,x_1)$ in contradition to Lemma~\ref{ab}.

Therefore $a \in VS(w_2,w_3) - \{w_2\}$. Then $b \in VS(w_2, a)$, for otherwise there would be a vertex of degree $2$ in $G$. 
First we consider the case that the trace of $S$ is $0$. We define 
$$ T'=T_1 \cup S(a,x_1) \cup A(x_1,b) \cup  T_2.$$
This is a directed path from $y_2$ to $x_2$ with 
$$T' \cap (S(y_1,a) - (A \cup B)) = \emptyset .$$
By Theorem~\ref{haupt} we have 
$$ G=G[A \cup B \cup S' \cup T'], $$
a contradiction, since 
$$  S(w_1,w_2)  \cap (A \cup B \cup S' \cup T') = \emptyset. $$

Now we consider the case that the trace of $S$ is $0'$. By Lemma~\ref{long}
and Lemma~\ref{next} the path $T_1(v,a)$ is an $\overline{A \cup B}$-arc.
We define the following directed minimal path $S''$ from $y_1$ to $x_1$:
$$ S''=A(y_1,v) \cup T_1(v,a) \cup S(a,x_1). $$
The trace of $S''$ is $0'$. Let $u$ be the origin of the 
$A(y_1,v_1)T_1$-arc with terminus $v$. Note that 
$$ T = T_1 \cup S(a,x_1) \cup A(x_1,y_1) \cup S(y_1,b) \cup T_2.$$ 
Then $u$ is the first 
vertex in $T$ that is also in $S''$ and $w_1$ is the last vertex in 
$T$ that is also in $S''$, since the trace of $T_2$ does not contain $2$.
If we replace $S$ by $S''$, this is a contradiction by Lemma~\ref{ab} and finally shows 
that the trace of $T_1$ does not contain $2$. 

Now we  know that the trace of $T_1$ is a string over $\{0,1\}$. By Lemma~\ref{long}
the trace of $T_1$ contains neither $00$ nor $11$, and by Lemma~\ref{next}
it contains neither $01$ nor $10$. Therefore the trace of $T_1$ is either 
empty, $0$ or $1$. Since the directed path $T_1$ has property $(A')$ the case that the trace of $T_1$ is 
$1$ is not possible. \qed

\begin{rem}
\label{ein}
An argument similar to the one that showed that the trace of $T_1$ does not contain 
$0'$ leads to the following observation:
If the trace of $S$ is $0'$, the symbol $0$ in the trace of $T_1$ corresponds 
to an $A(w_3,y_2)T_1$-arc and the symbol $0$ in the trace of $T_2$ corresponds to 
an $A(x_2,w_2)T_2$-arc.
\end{rem}

\begin{lem}
\label{b<a}
Suppose $b <_S a$. Then $S(b,a) \subseteq S(w_2,w_3)$. 
\end{lem}

{\it Proof}: The assertion can be deduced from Lemma~\ref{ab} 
after we show that $S(b,a) \subseteq A \cup B$. 

Define 
$$ X=S(y_1,b) \cup T_2$$
and 
$$ Y=T_1 \cup S(a,x_1).$$
Then $X$ is a directed path from $y_1$ to $x_2$ and $Y$ is a directed path 
from $y_2$ to $x_1$. Since $S(b,a) \cap (X \cup Y) = \emptyset$ and 
$$G=G[A \cup B \cup X \cup Y]$$
by Theorem~\ref{haupt}, we have $S(b,a) \subseteq A \cup B$. \qed

\begin{lem}
\label{sec}
If the trace of $T_1$ is empty then $a \notin VS(w_3,x_1)$.\\
If the trace of $T_1$ is $0$ then $a \notin VS(y_1,w_3)$.\\
If the trace of $T_2$ is empty then $b \notin VS(y_1,w_2)$.\\
If the trace of $T_2$ is $0$ then $b \notin VS(w_2,x_1)$.
\end{lem}

{\it Proof}: By symmetry it suffices to show the assertions for $T_1$.
The first assertion is an immediate consequence of Lemma~\ref{afternext}.

Now suppose that the trace of $T_1$ is $0$. If the trace of $S$ is $0$ then $a \notin VS(w_1,w_3)$  
by Lemma~\ref{afterlong}; if the trace of $S$ is $0'$ then $a \notin VS(w_1,w_3)$
by Lemma~\ref{long}. 

Suppose the trace of $S$ is $0$, and that $a \in VS(y_1,w_1)$. Let $z_1$, $z_2$, $z_3$ be the vertices of $VT_1 - \{y_2, a\}$
in the order in which they appear as $T_1$ is traced from $y_2$. Then 
  $$ G[(A \cup B \cup  S \cup  T_1) - (A(y_2, z_1) \cup  A(z_2, z_3) \cup  B(x_2, y_1))] $$
is an even subdivision of $K_{3,3}$, in contradiction to the 
fact that $G$ is minimal non-Pfaffian. If the trace of $S$ is $0'$ then $a \notin VS(y_1,w_1)$ by 
Lemma~\ref{after3} and Remark~\ref{ein}.\qed

\smallskip

The following lemma gives a complete list of graphs to be considered if $S \cap T \not= \emptyset$.

\begin{lem}
\label{notempty} 
Suppose $S \cap T \not= \emptyset$ and that the directed path $T$ has property $(A')$. Then one of the following cases is true:
\begin{enumerate}
\item the traces of $T_1$ and $T_2$ are $\Lambda$ and 0 respectively, $a \in VS(w_1,w_2) - \{w_2\}$, $b \in  VS(w_1,w_2) - \{w_1\}$, $a <_S b$,
\item the traces of $T_1$ and $T_2$ are empty, $a \in VS(y_1,w_1)$, $b \in VS(w_2,w_3)$,
\item the traces of $T_1$ and $T_2$ are empty, $a \in VS(y_1,w_1)$, $b \in VS(w_4,x_1)$, 
\item the traces of $T_1$ and $T_2$ are empty,  $a \in VS(w_2,w_3)$, $b \in VS(w_2,w_3)$, $a <_S b$, 
\item the traces of $T_1$ and $T_2$ are empty, $a \in VS(w_2,w_3)$, $b \in VS(w_2,w_3)$, $b <_S a$,
\item the traces of $T_1$ and $T_2$ are empty, $a \in VS(w_2,w_3)$, $b \in VS(w_4,x_1)$,
\item the traces of $T_1$ and $T_2$ are 0 and $\Lambda$ respectively, $a \in VS(w_3,w_4) - \{w_4\}$, $b \in  VS(w_3,w_4) - \{w_3\}$ and $a <_S b$.
\end{enumerate}
\end{lem}

{\it Proof}: First we deal with the case that
  $$ a \in (VS(w_1, w_2) - \{w_2\}) \cup (VS(w_3, w_4) - \{w_4\}) $$
or
  $$ b \in (VS(w_1, w_2) - \{w_1\}) \cup (VS(w_3, w_4) - \{w_3\}). $$  

Since there are no vertices of degree $2$ in $G$, 
the vertex $a$ is in $VS(w_1,w_2) - \{w_2\}$ if and only if 
$b \in VS(w_1,w_2) - \{w_1\}$. Furthermore $a <_S b$ in this case and by Lemma~\ref{sec} the trace of 
$T_1$ is empty and the trace of $T_2$ is $0$.
This situation corresponds to the first case in the lemma.

Similarly $a \in VS(w_3,w_4) - \{w_4\}$ if and only if 
$b \in VS(w_3,w_4) - \{w_3\}$. In this case $a <_S b$, the trace of 
    $T_1$ is $0$ and the trace of $T_2$ is empty.
This situation corresponds to the last case in the lemma. 

Therefore we may now assume that
$$ \{a,b\} \subset VS(y_1,w_1) \cup VS(w_2,w_3) \cup VS(w_4,x_1). $$

We show that the trace of $T_1$ is empty. Suppose the contrary, that the 
trace of $T_1$ is $0$. By Lemma~\ref{sec} and our assumption we have
$a \in VS(w_4,x_1)$. By Lemma~\ref{b<a} we have $a<_S b$ and therefore $b \in VS(w_4,x_1)$. This 
is a contradiction to Lemma~\ref{ab}. 

Similarly the trace of $T_2$ is empty.

By Lemma~\ref{sec} and our assumption  
we have
$$ a \in VS(y_1,w_1) \cup VS(w_2,w_3) $$
and 
$$b  \in VS(w_2,w_3) \cup VS(w_4,x_1). $$

From this result we deduce the following list of cases to be considered.
\begin{enumerate}
\item $a \in VS(y_1,w_1)$, $b \in VS(w_2,w_3)$
\item $a \in VS(y_1,w_1)$, $b \in VS(w_4,x_1)$ 
\item $a \in VS (w_2,w_3)$, $b \in VS(w_2,w_3)$, $a <_S b$
\item $a \in VS (w_2,w_3)$, $b \in VS(w_2,w_3)$, $b <_S a$
\item $a \in VS(w_2,w_3)$, $b \in VS(w_4,x_1)$.
\end{enumerate}
\qed

\begin{rem}
\label{symmetric2}
If we change from the reference orientation with respect to $(M,N,f)$ to the reference orientation with respect to 
$(N,M,f)$ the first case in Lemma~\ref{notempty} changes to the last case. Therefore we do not have to consider the 
last case. The same is true for the second case and the sixth case in Lemma~\ref{notempty}, and so we do not consider the sixth case. 
\end{rem}

\section{The minimal non-Pfaffian near bipartite graphs}

In this section we consider the cases in Lemma~\ref{nointersection} and in Lemma~\ref{notempty}, and with this
complete the proof of Theorem~\ref{main}. For that purpose we need the following lemma which has already been 
proved in \cite{LiReFi}.

\begin{lem}
\label{contract}
  Let $G$ be a graph with a circuit $C$ of odd length and let $G^C$ be the
  graph obtained from $G$ by contracting $VC$. If $G^C$ is not Pfaffian,
  then neither is $G$.
\end{lem}

We divide the argument into cases according to whether or not $S \cap T = \emptyset$.

{\it Case 1}: In the case where $S \cap T = \emptyset$ it follows from Lemma~\ref{nointersection} and Remark~\ref{symmetric1} that the trace of $S$ may be assumed to be  
either $0$ or $0'$
and that of $T$ may be assumed to be one of $0$, $0'$, $20$, $21$. Let the vertices in $VT - \{y_2, x_2\}$ be $z_1, z_2, \dots, z_n$ in the order in which they appear. 

{\parindent0cm \it Subcase 1.1}: Suppose the trace of $T$ is 0.

{\parindent0cm \it Subcase 1.1.1}: Suppose the trace of $S$ is 0. Consider the circuits
\begin{eqnarray*}
C_1 &=& S(w_1, w_2) \cup A(w_2, w_1), \\
C_2 &=& S(w_3, w_4) \cup A(w_4, w_3), \\
C_3 &=& T(z_1, z_2) \cup A(z_2, z_1), \\
C_4 &=& T(z_3, z_4) \cup A(z_4, z_3), \\
C_5 &=& S \cup B(x_1, y_2) \cup T \cup B(x_2, y_1).
\end{eqnarray*}
Their sum is $A + B$. However, under our reference orientation all of $C_1$, $C_2$, $C_3$, $C_4$, $C_5$, $A + B$ are clockwise even, 
but under our extended Pfaffian orientation
of $H$ only $A + B$ is clockwise even, by Corollary~\ref{A+B}. This result contradicts Lemma~\ref{parity}.
{\parindent0cm \it Subcase 1.1.2}: Suppose the trace of $S$ is $0'$. By symmetry we may assume that $w_3 <_{A(x_2, y_2)} z_2$. 
Then $G$ is isomorphic to $\Gamma_1$. The 
isomorphism $\phi$ from $\Gamma_1$ and $G$ is given by
$\phi(a) = x_2$, $\phi(b) = y_1$, $\phi(c) = w_1$, $\phi(d) = w_2$, $\phi(e) = w_3$, $\phi(f) = w_4$, $\phi(g) = x_1$, $\phi(h) = y_2$, $\phi(i) = z_1$, $\phi(j) = z_2$,
$\phi(k) = z_3$, $\phi(l) = z_4$.

\medskip

{\parindent0cm \it Subcase 1.2}: Suppose the trace of $T$ is $0'$.

{\parindent0cm \it Subcase 1.2.1}: The case where the trace of $S$ is 0 is symmetric to Subcase~1.1.2.

{\parindent0cm \it Subcase 1.2.2}: Suppose the trace of $S$ is $0'$. Consider the circuits
\begin{eqnarray*}
C_1 &=& S \cup A(x_1, y_1), \\
C_2 &=& T \cup A(x_2, y_2), \\
C_3 &=& S \cup B(x_1, y_2) \cup T \cup B(x_2, y_1). 
\end{eqnarray*}
We obtain a contradiction by the method in Subcase~1.1.1.

\medskip

{\parindent0cm \it Subcase 1.3}: Suppose the trace of $T$ is 20. 

{\parindent0cm \it Subcase 1.3.1}: Suppose the trace of $S$ is 0. This case is similar to Subcase~1.1.1 except that vertices $z_2$, $z_3$, $z_4$ in Subcase~1.1.1 are replaced
by $z_4$, $z_5$, $z_6$ respectively.

{\parindent0cm \it Subcase 1.3.2}: Suppose the trace of $S$ is $0'$. Note that $z_1 = v_2$ and $z_3 = v_1$. Contract the circuit
$T(z_1, z_3) \cup \{e_2\}$. The resulting graph is isomorphic to $\Gamma_1$ whether or not $w_3 <_{A(x_2, y_2)} z_4$. 

\medskip

{\parindent0cm \it Subcase 1.4}: Suppose the trace of $T$ is $21$.

{\parindent0cm \it Subcase 1.4.1}: Suppose the trace of $S$ is 0. Consider the circuits
\begin{eqnarray*}
C_1 &=& S(w_1, w_2) \cup A(w_2, w_1), \\
C_2 &=& S(w_3, w_4) \cup A(w_4, w_3), \\
C_3 &=&  T \cup A(x_2, y_2), \\
C_4 &=& T(y_2, z_5) \cup B(z_5, y_2), \\
C_5 &=& S \cup B(x_1, z_5) \cup T(z_5, x_2) \cup B(x_2, y_1). 
\end{eqnarray*}
We obtain a contradiction by the method in Subcase~1.1.1.

{\parindent 0cm \it Subcase 1.4.2}: Suppose the trace of $S$ is $0'$.  Contract the circuit
$T(z_1, z_3) \cup \{e_2\}$. The resulting graph is isomorphic to $\Gamma_1$.

\medskip

{\it Case 2}: If $S \cap T \neq \emptyset$ then by Remark~\ref{symmetric2} we see that only cases 1--5 in Lemma~\ref{notempty} need to be considered.
Case $i$ of the lemma is dealt with in Subcase 2.$i$ below. Let the vertices in $(VT_1 \cup VT_2) - \{y_2, x_2, a, b\}$
be $z_1, z_2, \dots, z_n$ in the order in which they appear in $T$. 

{\parindent0cm \it Subcase 2.1}: The traces of $T_1$ and $T_2$ are $\Lambda$ and 0 respectively, $a \in VS(w_1,w_2) - \{w_2\}$, $b \in  VS(w_1,w_2) - \{w_1\}$, $a <_S b$.

{\parindent0cm \it Subcase 2.1.1}: Suppose the trace of $S$ is 0. Note that $z_1 = v_2$ and $w_1 = v_1$. Contraction of the circuit 
$T(z_1, a) \cup S(w_1, a) \cup \{e_2\}$ yields a graph isomorphic to $\Gamma_1$. 

{\parindent0cm \it Subcase 2.1.2}: Suppose the trace of $S$ is $0'$. By Remark~\ref{ein} we have $z_3 <_{A(x_2, y_2)} w_2$. Define 
  $$ S' = S(y_1, b) \cup T(b, z_3) \cup A(z_3, w_2) \cup S(w_2, x_1). $$
By Theorem~\ref{haupt} we find that $G = G[A \cup B \cup S' \cup T]$, in contradiction to the fact that 
  $$ S(b, w_2) \cap (A \cup B \cup S' \cup T) = \emptyset. $$

{\parindent0cm \it Subcase 2.2}: The traces of $T_1$ and $T_2$ are empty, $a \in VS(y_1,w_1)$, $b \in VS(w_2,w_3)$. Note that $a = y_1$ and $b = w_3$.

{\parindent0cm \it Subcase 2.2.1}: Suppose the trace of $S$ is 0. Consider the circuits 
\begin{eqnarray*}
C_1 &=& S(w_1, w_2) \cup A(w_2, w_1), \\
C_2 &=& S(w_3, w_4) \cup A(w_4, w_3), \\
C_3 & = & T_1 \cup S(y_1, w_3) \cup T_2 \cup A(x_2, y_2), \\
C_4 & = & T_1 \cup S \cup B(x_1, y_2), \\
C_5 & = & S(y_1, w_3) \cup T_2 \cup B(x_2, y_1).
\end{eqnarray*}
We obtain a contradiction by the method in Subcase~1.1.1.       

{\parindent0cm \it Subcase 2.2.2}: Suppose the trace of $S$ is $0'$. Consider the circuits       
\begin{eqnarray*}
C_1 &=& S \cup A(x_1, y_1), \\
C_2 & = & T_1 \cup S(y_1, w_3) \cup A(w_3, y_2), \\
C_3 & = & T_2 \cup A(x_2, w_3), \\
C_4 & = & T_1 \cup S \cup B(x_1, y_2), \\
C_5 & = & S(y_1, w_3) \cup T_2 \cup B(x_2, y_1).
\end{eqnarray*}
We obtain a contradiction by the method in Subcase~1.1.1.

{\parindent0cm \it Subcase 2.3}: The traces of $T_1$ and $T_2$ are empty, $a \in VS(y_1,w_1)$, $b \in VS(w_4,x_1)$. Note that $a = y_1$ and $b = x_1$.

{\parindent0cm \it Subcase 2.3.1}: Suppose the trace of $S$ is 0. Consider the circuits 
\begin{eqnarray*}
C_1 &=& S(w_1, w_2) \cup A(w_2, w_1), \\
C_2 &=& S(w_3, w_4) \cup A(w_4, w_3), \\
C_3 & = & T_1 \cup S \cup T_2 \cup A(x_2, y_2), \\
C_4 & = & T_1 \cup S \cup B(x_1, y_2), \\
C_5 & = & S \cup T_2 \cup B(x_2, y_1).
\end{eqnarray*}
We obtain a contradiction by the method in Subcase~1.1.1.

{\parindent0cm \it Subcase 2.3.2}: Suppose the trace of $S$ is $0'$. Consider the circuits
\begin{eqnarray*}
C_1 &=& S \cup A(x_1, y_1), \\
C_2 & = & T_1 \cup S(y_1, w_3) \cup A(w_3, y_2), \\
C_3 & = & S(w_2, x_1) \cup T_2 \cup A(x_2, w_2), \\
C_4 & = & T_1 \cup S \cup B(x_1, y_2), \\
C_5 & = & S \cup T_2 \cup B(x_2, y_1).
\end{eqnarray*}
We obtain a contradiction by the method in Subcase~1.1.1.

{\parindent0cm \it Subcase 2.4}: The traces of $T_1$ and $T_2$ are empty, $a \in VS(w_2,w_3)$, $b \in VS(w_2,w_3)$, $a <_S b$. Note that $a = w_2$ and $b = w_3$.

{\parindent0cm \it Subcase 2.4.1}: Suppose the trace of $S$ is 0. Consider the circuits 
\begin{eqnarray*}
C_1 &=& S(w_1, w_2) \cup A(w_2, w_1), \\
C_2 &=& S(w_3, w_4) \cup A(w_4, w_3), \\
C_3 & = & T_1 \cup S(w_2, w_3) \cup T_2 \cup A(x_2, y_2), \\
C_4 & = & T_1 \cup S(w_2, x_1) \cup B(x_1, y_2), \\
C_5 & = & S(y_1, w_3) \cup T_2 \cup B(x_2, y_1).
\end{eqnarray*}
We obtain a contradiction by the method in Subcase~1.1.1.

{\parindent0cm \it Subcase 2.4.2}: Suppose the trace of $S$ is $0'$. This case is symmetric to Subcase 2.4.1.

{\parindent0cm \it Subcase 2.5}: The traces of $T_1$ and $T_2$ are empty, $a \in VS(w_2,w_3)$, $b \in VS(w_2,w_3)$, $b <_S a$.

{\parindent0cm \it Subcase 2.5.1}: Suppose the trace of $S$ is 0. Then $G$ is isomorphic to $\Gamma_2$.

{\parindent0cm \it Subcase 2.5.2}: Suppose the trace of $S$ is $0'$. Consider the circuits 
\begin{eqnarray*}
C_1 & = & S \cup A(x_1, y_1), \\
C_2 & = & T_1 \cup A(a, y_2), \\
C_3 & = & T_2 \cup A(x_2, b), \\
C_4 & = & T_1 \cup S(a, x_1) \cup B(x_1, y_2), \\
C_5 & = & S(y_1, b) \cup T_2 \cup B(x_2, y_1).
\end{eqnarray*}
We obtain a contradiction by the method in Subcase~1.1.1.
The proof of Theorem~\ref{main} is now complete.

\section{Non-reduction of $\Gamma_1$ and $\Gamma_2$ to $K_{3,3}$}

We conclude the paper by showing that neither $\Gamma_1$ nor $\Gamma_2$ is reducible to an even subdivision of $K_{3,3}$.

\begin{lem}
Both $\Gamma_1$ and $\Gamma_2$ are minimal non-Pfaffian graphs.
\end{lem}

\begin{figure}
\begin{center}
\setlength{\unitlength}{1cm}
\begin{picture}(14.5,7)
\put(0,0)
{\scalebox{0.35}{\includegraphics{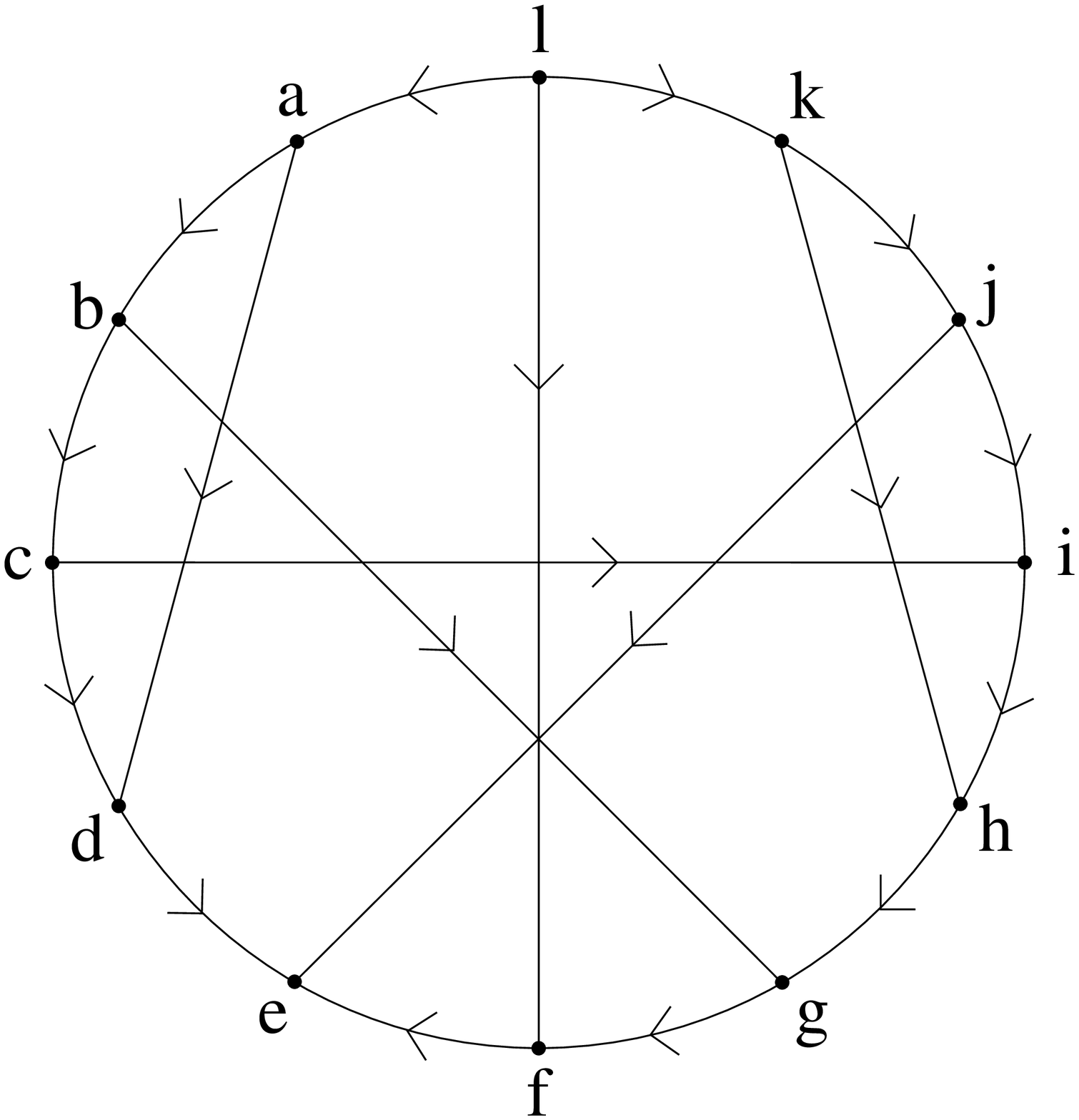}}}
\put(8,0)
{\scalebox{0.35}{\includegraphics{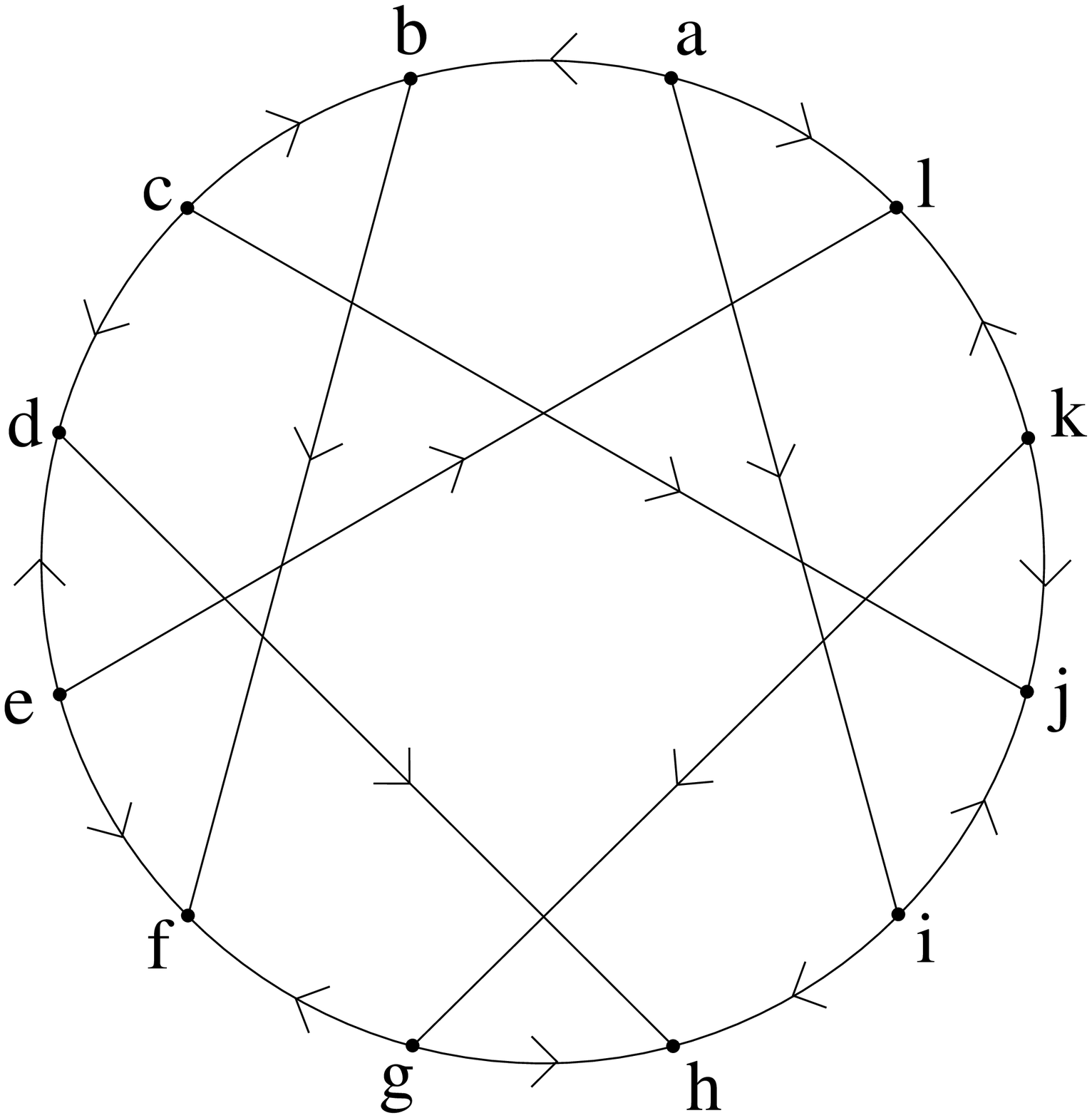}}}
\makebox(7,-1){$\Gamma_1$}
\makebox(8.5,-1){$\Gamma_2$}
\end{picture}
\end{center}
\caption{}
\label{reorient}
\end{figure}

{\it Proof}: Let $\Gamma_1$ and $\Gamma_2$ be oriented as in Figure~\ref{reorient}.

First we consider $\Gamma_1$, which we have already seen to be non-Pfaffian. Suppose therefore that $\Gamma_1$ is not minimal.
Let $x$ be an edge such that $\Gamma_1 - \{x\}$ is non-Pfaffian.
The $1$-factors of $\Gamma_1$ are 
\begin{eqnarray*}
f_1 & = & \{(a,b),(c,d),(f,e),(h,g),(j,i),(l,k)\}, \\
f_2 & = & \{(l,a),(b,c),(d,e),(g,f),(i,h),(k,j)\}, \\
f_3 & = & \{(b,c),(a,d),(f,e),(h,g),(j,i),(l,k)\}, \\
f_4 & = & \{(a,b),(c,d),(j,e),(g,f),(i,h),(l,k)\}, \\
f_5 & = & \{(a,d),(b,c),(j,e),(g,f),(i,h),(l,k)\}, \\
f_6 & = & \{(l,a),(b,c),(d,e),(g,f),(k,h),(j,i)\}, \\
f_7 & = & \{(l,a),(b,g),(i,h),(k,j),(c,d),(f,e)\}, \\
f_8 & = & \{(l,a),(b,g),(k,h),(c,d),(f,e),(j,i)\}, \\
f_9 & = & \{(l,f),(d,e),(c,i),(k,j),(a,b),(h,g)\}, \\
f_{10} & = & \{(a,d),(b,g),(c,i),(j,e),(k,h),(l,f)\}. 
\end{eqnarray*}
Observe that the figure for the undirected graph $\Gamma_1$ is symmetric about 
the edge $(l,f)$. Therefore we can assume that $x \notin f_1$. All 1-factors 
are associated with a plus sign except $f_{10}$. Thus $x \notin f_{10}$, for 
otherwise $\Gamma_1 - \{x\}$ is Pfaffian, and therefore 
$x \in f_2$. 

Suppose that $x=(d,e)$ or $x=(k,j)$. We obtain a Pfaffian orientation 
of $\Gamma_1 - \{x\}$ if we change the orientation of $(l,f)$, since 
every 1-factor of $G$ that contains $(l,f)$ also contains $x$ except 
for $f_{10}$. If $x=(g,f)$ or $x=(i,h)$,  we obtain a Pfaffian orientation 
of $\Gamma_1 - \{x\}$ by changing the orientation of $(j,e)$. 
If $x=(l,a)$,  we obtain a Pfaffian orientation 
of $\Gamma_1 - \{x\}$ by changing the orientation of $(b,g)$. 
If $x=(b,c)$,  we obtain a Pfaffian orientation of 
$\Gamma_1 - \{x\}$ by changing the orientation of $(a,d)$.
In all cases we have a contradiction to the fact that $\Gamma_1 - \{x\}$
was non-Pfaffian. Therefore $\Gamma_1$ is minimal non-Pfaffian.

\smallskip

Now suppose that the non-Pfaffian graph $\Gamma_2$ is not minimal, and let $x$ be 
an edge such that  $\Gamma_2 - \{x\}$ is non-Pfaffian. The $1$-factors of 
$\Gamma_2$ are
\begin{eqnarray*}
f_1 & = & \{(a,b),(c,d),(e,f),(g,h),(i,j),(k,l)\}, \\
f_2 & = & \{(c,b),(e,d),(g,f),(i,h),(k,j),(a,l)\}, \\
f_3 & = & \{(a,b),(c,j),(e,d),(g,f),(i,h),(k,l)\}, \\
f_4 & = & \{(a,b),(c,d),(e,l),(g,f),(i,h),(k,j)\}, \\
f_5 & = & \{(c,b),(d,h),(k,g),(a,l),(e,f),(i,j)\}, \\
f_6 & = & \{(b,f),(e,d),(c,j),(a,i),(k,l),(g,h)\}, \\
f_7 & = & \{(b,f),(a,i),(k,j),(e,l),(c,d),(g,h)\}, \\
f_8 & = & \{(e,d),(b,f),(a,l),(k,g),(i,h),(c,j)\}, \\
f_9 & = & \{(g,f),(d,h),(c,b),(a,i),(k,j),(e,l)\}, \\
f_{10} & = & \{(b,f),(e,l),(k,g),(d,h),(c,j),(a,i)\}.
\end{eqnarray*}
Observe that there is an automorphism of $\Gamma_2$ that interchanges the $f_1$-alternating circuits 
$f_1+f_2$ and $f_1 + f_{10}$. Therefore we can assume that $x \in f_1 + f_{10}$. The figure 
for the undirected graph $\Gamma_2$ is symmetric about the line through the midpoints of 
the edges $(a,b)$ and $(g,h)$. Therefore we can assume that 
\begin{eqnarray*}
x \in (f_1 + f_{10}) - \{(b,f),(k,g),(c,j)\} =  \\
\{(a,b),(c,d),(e,f),(g,h),(i,j),(k,l),(e,l),(d,h),(a,i)\}.
\end{eqnarray*}
If $x \in f_1$ then the given orientation is a Pfaffian orientation of 
$\Gamma_2 - \{x\}$, since the sign of $f_1$ is the opposite of that of the other 1-factors.
Therefore $x \in \{(e,l),(d,h),(a,i)\}$. If $x=(e,l)$, we obtain a Pfaffian orientation of 
$\Gamma_2 - \{x\}$ by changing the orientation of $(c,d)$, since every $1$-factor that contains $(c,d)$ also contains 
$x$ except for $f_1$. If $x=(d,h)$, we obtain a Pfaffian orientation 
of $\Gamma_2 - \{x\}$ by changing the orientation of $(e,f)$. If $x=(a,i)$, we obtain a Pfaffian orientation 
of $\Gamma_2 - \{x\}$ by changing the orientation of $(g,h)$. In all cases we have a contradiction to the 
fact that $\Gamma_2 - \{x\}$ was Pfaffian.  Therefore $\Gamma_2$ is minimal non-Pfaffian.

\begin{cor}
Neither $\Gamma_1$ nor $\Gamma_2$ contains an even subdivision of $K_{3,3}$.
\end{cor}

{\it Proof}: If $\Gamma_1$ or $\Gamma_2$ contained an even subdivision of $K_{3,3}$ then $\Gamma_1$ or $\Gamma_2$ 
itself would be an even subdivision of 
$K_{3,3}$, since $\Gamma_1$ and $\Gamma_2$ are minimal non-Pfaffian and every even subdivision of $K_{3,3}$ is non-Pfaffian.
But $\Gamma_1$ and $\Gamma_2$ both have 12 vertices of degree 3 whereas an even subdivision of $K_{3,3}$ has only 6.  \qed  

In order to see that $\Gamma_1$ and $\Gamma_2$ are not reducible to 
an even subdivision of $K_{3,3}$, we need the following two lemmas.

\begin{lem}
\label{formel}
Let $G$ be a minimal non-Pfaffian graph and $C$ a circuit of odd length in $G$. 
Suppose that the graph $G^C$ obtained by contracting $VC$ to a vertex $v$ is also non-Pfaffian.
Let $VC=\{v_1,v_2,\dots,v_n\}$. 
Then
$$ \deg v = \sum_{i=1}^{n} \deg v_i -2 n. $$
Moreover if $w$ is a vertex in $VG - VC$ then $w$ is also a vertex of $G^C$, and 
$$ \deg_{G} w = \deg_{G^C} w. $$
\end{lem}

{\it Proof}:
Let $N(u)$ denote the set of vertices in $G$ that are adjacent to the vertex $u$ in $G$.
Moreover we assume that $v_i$ is adjacent to $v_{i+1}$ in $C$, for all $i < n$. We define $v_0=v_n$
and $v_{n+1}=v_1$. 

First we observe the following. Let $i,j \in \{1,2,\dots,n\}$ and $i \not= j$. We claim that   
$$ N(v_i)  \cap N(v_j)  \subseteq \{v_{i-1},v_{i+1}\}.$$
Indeed, suppose $u \in N(v_i) \cap N(v_j)$. Then there are edges $e_i$ and $e_j$ joining $u$ to $v_i$ and $v_j$ respectively. If $u \notin VC$, then 
the graph obtained from $G - \{e_i\}$ by contracting $VC$ is $G^C$. 
Since $G^C$ is supposed to be non-Pfaffian, it follows from Lemma~\ref{contract}
that $G - \{e_i\}$ is non-Pfaffian too. This is a contradiction to the 
fact that $G$ was minimal non-Pfaffian. Therefore 
$u \in VC$. If $e_i \notin C$ or $e_j \notin C$ we can 
conclude in a similar way that either  $G - \{e_i\}$ 
or $G - \{e_j\}$ is non-Pfaffian, and therefore have a contradiction again.
Thus $\{e_i, e_j\} \subseteq C$, and the claim is proved. We infer that  $C$ has no chords, and any vertex not in $C$ has at most one neighbour in $C$. 
We have
$$ N(v) = \bigcup_{i=1}^n N(v_i) - VC = \bigcup_{i=1}^n (N(v_i) - \{v_{i-1},v_{i+1}\}),$$
and so 
$$ \deg v = |N(v)| = \sum_{i=1}^n |N(v_i) - \{v_{i-1},v_{i+1}\}| = \sum_{i=1}^n (|N(v_i)|-2) = 
   \sum_{i=1}^{n} \deg v_i -2 n. $$

Finally the degree of a vertex $w$ that is not in $C$ does not change upon contraction of $C$, since 
$w$ is adjacent to at most one vertex in $C$. \qed

\begin{lem}
\label{last}
Let $G$ be a minimal non-Pfaffian graph that is cubic and does not contain a circuit of length 
$3$. Then $G$ is not reducible to an even subdivision of $K_{3,3}$.
\end{lem}

{\it Proof}: Suppose the contrary, that is that there exists a sequence 
$C_0, C_1, \dots, C_{p-1}$ of circuits of odd length and a sequence $G_0, G_1, \dots, G_p$ of graphs such
that $G_0=G$, $G_p$ is an even subdivision of $K_{3,3}$ and, for all $i < p$, $G_{i+1}$ is obtained from 
$G_i$ by contracting $VC_i$. We see inductively that  for all $i$ the graph $G_i$ is minimal non-Pfaffian
by Lemma~\ref{contract}, since $G_p$ is non-Pfaffian and $G_0$ is minimal non-Pfaffian. 

First we show that $G_1$ contains a vertex of degree at least $5$. Since $G_0$ does not 
contain a circuit of length $3$, the circuit $C_0$ must have length at least 5. 
Let $VC_0=\{v_1,v_2,\dots,v_n\}$ and let $v$ be the corresponding vertex in $G_1$. 
Then, by Lemma~\ref{formel}, 
$$ \deg v = \sum_{i=1}^n \deg v_i - 2 n = 3 n - 2 n = n \ge 5. $$
Again by Lemma~\ref{formel} all the other vertices in $G_1$ have degree $3$, since 
$G_0$ is cubic. 

Now we show by induction that all the vertices in $G_i$, where $i \ge 1$,  
are of degree at least $3$ and that there exists a vertex in $G_i$ with 
degree at least $5$. 
Therefore let us assume that the induction hypothesis is true for $G_i$ and 
show it for $G_{i+1}$. First we show that every vertex in $G_{i+1}$ is  of degree at least $3$. For all vertices in $G_{i+1}$ except the one that corresponds to $C_i$
this is an immediate consequence of Lemma~\ref{formel} and the induction hypothesis. 
Let $w$ be the vertex in $G_{i+1}$ that corresponds to $C_i$ and let
$VC_i=\{w_1,w_2,\dots,w_m\}$. Then, by Lemma~\ref{formel} and the induction hypothesis,
$$ \deg w = \sum_{i=1}^m \deg w_i -2 m  \ge 3 m - 2 m =  m \ge 3.$$
Now let $u$ be the vertex of degree at least $5$ in $G_i$. If $u \notin VC_i$
then $u \in VG_{i+1}$ and $\deg_{G_{i+1}} u = \deg_{G_i} u $, by Lemma~\ref{formel}. 
Therefore suppose that $u \in VC_i$, and without loss of generality assume 
that $u=w_1$. Then, by Lemma~\ref{formel} and the induction hypothesis,
$$ \deg w = \deg u + \sum_{i=2}^m \deg w_i - 2 m \ge \deg u + 3 (m-1) - 2 m = 
\deg u + m - 3 \ge \deg u. $$

Therefore $G_p$ contains a vertex of degree at least $5$. This is a contradiction, since 
$G_p$ was an even subdivision of $K_{3,3}$. \qed 

\begin{cor}
Neither $\Gamma_1$ nor $\Gamma_2$ is reducible to an even subdivision of $K_{3,3}$.
\end{cor}

{\it Proof}: This is an immediate consequence of Lemma~\ref{last} since both $\Gamma_1$ 
and $\Gamma_2$ are minimal non-Pfaffian, cubic and do not contain a circuit of 
length $3$. \qed

\end{document}